\documentclass[german,english,reqno]{amsart}
\usepackage{bookman}
\usepackage[T1]{fontenc}
\usepackage[latin1]{inputenc}
\usepackage{variorefinc}
\usepackage{amssymb}
\IfFileExists{url.sty}{\usepackage{url}}
                      {\newcommand{\url}{\texttt}}

\makeatletter

\providecommand{\LyX}{L\kern-.1667em\lower.25em\hbox{Y}\kern-.125emX\@}

 \theoremstyle{definition}
 \newtheorem*{defn*}{Definition}
 \theoremstyle{definition}
  \newtheorem*{example*}{Example}
 \theoremstyle{plain}    
 \newtheorem{lem}{Lemma} 
 \theoremstyle{remark}
 \newtheorem*{rem*}{Remark}
 \theoremstyle{plain}    
 \newtheorem*{thm*}{Theorem} 

\usepackage{color}
\usepackage[marginal]{footmisc}
\usepackage{graphicx}
\usepackage{graphicx}

\newcommand{\ZZ}{\mathbb Z}
\newcommand{\NN}{\mathbb N}
\newcommand{\RR}{\mathbb R}
\newcommand{\CC}{\mathbb C}
\newcommand{\QQ}{\mathbb Q}
\newcommand{\DD}{\mathbb D}
\newcommand{\TT}{\mathbb T}

\DeclareMathOperator{\EE}{\mathbb E}
\DeclareMathOperator{\PP}{\mathbb P}
\DeclareMathOperator{\real}{Re}
\DeclareMathOperator{\imag}{Im}
\DeclareMathOperator{\diam}{diam}

\newcommand{\one}{\mathbf {1}}

\newcommand{\half}{{\textstyle \frac{1}{2}}}
\newcommand{\third}{{\textstyle \frac{1}{3}}}

\newcommand{\quarter}{{\textstyle \frac{1}{4}}}

\newcommand{\brk}{\\ }

\newcommand{\sh}{ }

\newcommand{\enspl}{\end{split}}
\def\tsh#1{{\textstyle #1}}

\newcommand{\Meas}{\mathfrak {M}}
\newcommand{\Haus}{\mathfrak {H}}

\DeclareMathOperator{\LE}{LE}
\DeclareMathOperator{\QL}{QL}
\DeclareMathOperator{\str}{str}
\DeclareMathOperator{\supp}{supp}

\newcommand{\dreg}{\mathcal{D}}
\newcommand{\ereg}{\mathcal{E}}

\newcommand{\eclass}{\mathcal{X}}
\newcommand{\ieeclass}{\mathcal{Y}}

\newcommand{\brmul}{\discretionary{\mbox{$\,\cdot$}}{}{}}

\allowdisplaybreaks

\usepackage{babel}
\makeatother
\begin{document}

\title{Scaling limit of loop erased random walk --- a naive approach}

\author{Gady Kozma}

\thanks{This work is part of the research program of the European Network
{}``Analysis and Operators'', contract HPRN-CT-00116-2000 supported
by the European Commission.}

\curraddr{The Weizmann Institute of Science, Rehovot, Israel.}

\email{gadykozma@hotmail.com, gadyk@wisdom.weizmann.ac.il}

\begin{abstract}
We give an alternative proof of the existence of the scaling limit
of loop-erased random walk which does not use L\"owner's differential
equation.
\end{abstract}
\maketitle
\theoremstyle{remark}
\newtheorem*{rmks}{Remarks}
\theoremstyle{plain}
\newtheorem{sublem}{Sublemma}[lem] 
\newtheorem*{mainlem}{Main lemma}
\renewcommand{\labelenumi}{\theenumi .}
\newcounter{const}
\newcounter{Const}
\def\newc#1{
\refstepcounter{const}
\label{#1}
}
\def\newC#1{
\refstepcounter{Const}
\label{#1}
}

\section{Introduction}

Loop erased random walk is a process for creating a random simple
path, which starts from a regular random walk and then removes all
loops in a chronological order until a simple path is reached. In
dimension 2, it is typical to stop the process on the boundary of
some bounded domain $\dreg $, so the process creates a random simple
path from the point of origin to $\partial \dreg $. Originally \cite{L80}
it was suggested as a model for investigating the self-avoiding random
walk (i.e.~a random walk conditioned not to hit itself) but it was
found that these processes are cosingular. Notwithstanding, loop-erased
random walk is still a useful model for a random simple path. See
\cite{D92} for connections with various physical models such as the
{}``$Q$-states Potts model''%
\footnote{Loop-erased random walk is related to the case $Q=0$. It might be
interesting to note that critical percolation is also a particular
case, when $Q=1$.%
} and polymer coalescence. Another connection to physics which is also
interesting mathematically is the {}``Laplacian random walk,'' defined
in \cite{LEP86} and proved in \cite{L87} to be identical to loop-erased
random walk. The connection between loop-erased random walk and the
{}``uniform random spanning tree'' --- the spanning tree of a graph
chosen among all spanning trees with equal probabilities --- has given
thrust to the research of both. See \cite{P91,W96}.%
\footnote{The strongest result in this direction, Wilson's algorithm \cite{W96},
is stated in lemma \ref{lem_wilson} below.%
} The introduction to \cite{S00} explains all these connections in
a clear and concise way.

It is natural to assume that the distributions of loop-erased random
walks on the graphs $\dreg \cap \delta \ZZ ^{2}$ converge to a scaling
limit as $\delta \rightarrow 0$ which would be a {}``loop-erased
Brownian motion'' though this term per se is meaningless as the process
of loop erasure cannot be applied to Brownian motion: it has a dense
set of loops which cannot be ordered chronologically. Like many similar
processes, and in particular because regular random walk exhibits
this phenomenon, one might expect the limit to be conformally invariant.
As a rule of the thumb, conformal invariance can be expected for any
process which is local and invariant to scaling and rotation, since
a conformal map is, infinitesimally, just that, a rotation and scaling.%
\footnote{This might be the place to remark that loop-erased random walk is
formally \textbf{not} a local process, which is a major obstacle to
its analysis.%
} This conjecture lay open for a long period, with the first important
step done by Richard Kenyon \cite{K00,K00b} who proved the conformal
invariance of certain measurables of loop-erased walk, as well as
calculating explicit growth exponentials. Oded Schramm \cite{S00}
demonstrated how to describe the scaling limit of loop-erased random
walk using L\"owner's differential equation, assuming that the limit
exists and is conformally invariant. Basically he showed that the
generating function of L\"owner's equation is distributed like $e^{iB(2t)}$
where $B$ is a one dimensional Brownian motion (a good source%
\footnote{Notwithstanding the fact that Ahlfors' use of L\"owner's method for
the proof of Biberbach's conjecture for the case $n=3$ is a little
outdated.%
} on L\"owner's equation is \cite{A73}). This result opened the road
to the first proof of the conjecture \cite{LSW02}, and to additional
exciting results that connect other random processes to SLE (stochastic
L\"owner equation) with only a different multiplicative parameter
--- see \cite{S01,LSW02} for details. The aim of this paper is to
give an alternative proof of the existence of loop-erased random walk.

Why give another proof of a known result, and a longer one to boot?
Lawler-Schramm-Werner's proof is of the kind that {}``knows the answer''.
Very roughly, they started from the generating function of L\"owner's
equation for the discrete process (i.e.~the loop-erased random walk,
considered as a path in $\CC $ from $0$ to $\partial \dreg $),
showed that its distribution converges to Brownian motion as $\delta \rightarrow 0$
and then used compactness arguments to get convergence in the stronger
topology of simple paths in $\dreg $. My technique is {}``naive'',
it shows that loop-erased random walk converges without proving anything
about the limit. Thus, for example, it does not really distinguish
between simply connected and finitely connected domains.%
\footnote{For infinitely connected domains other factors are at work and the
loop-erased random walk does not necessarily converge to a limit.
See example on page \pageref{page:exam_punct}.%
} Each approach can be extended in directions the other cannot. At
the end of chapter \ref{chap_limit} we discuss very briefly and without
proofs some directions where this approach can be carried to.

I wish to thank Oded Schramm for reading early versions of this paper.

\subsection{About the proof}

Despite its length, in essence it is a simple proof, with the core
argument being localization and symmetry. Let $R$ be a random walk
on $\ZZ ^{2}$ from $0$ stopped on $\partial \dreg $ for some domain
$\dreg $. Let $S\subset \dreg $ be some (small) square. We write
\[
\LE (R)=\gamma _{1}\cup \gamma _{2}\cup \gamma _{3}\]
 where $\gamma _{1}$ is the portion of $\LE (R)$ until the first
time when $\LE (R)$ hits $S$. Notice that this is \textbf{not} the
same as the loop-erasure of a random walk stopped on $\partial S$!
$\gamma _{2}$ is the portion of $\LE (R)$ until the last time when
$\LE (R)$ is inside $S$, and $\gamma _{3}$ is the reminder (the
precise form of this division is in the main lemma, page \pageref{page:gam123}).
Tracing the process of loop-erasure in $\dreg $ one sees that $\gamma _{1}$
does not depend on anything that happens inside $S$: when one knows
all entry and exit points of $R$ from $S$, and all the trajectories
that $R$ does outside $S$, one can calculate $\gamma _{1}$. In
particular, if we compare random walks $R_{1}$ and $R_{2}$ on graphs
$G_{1}$ and $G_{2}$, where $G_{1}\setminus S=G_{2}\setminus S$
and inside $S$ we have some estimate of the sort\begin{equation}
p_{1}(v)\simeq p_{2}(v)\label{eq:pisketch}\end{equation}
where $p_{i}(v)$ is the probability of a random walk on $G_{i}$
to exit $S$ in a particular vertex $v$, then we should have that
\[
\gamma _{1,1}\simeq \gamma _{1,2}\quad .\]
This argument and the precise meaning of {}``$\simeq $'' are contained
in lemma \ref{lem_clerw_local}. To make this argument work for $\gamma _{3}$,
we have to use the symmetry of loop-erased random walk (exact details
in the main lemma). $\gamma _{2}$ describes what was coined in \cite{S00}
a {}``quasi-loop,'' and can be estimated using the methods ibid.
(see lemma \ref{lem_no_quasi_loops}).

This concludes the main argument, and leaves us with the question:
what are those mysterious graphs $G_{i}$ which differ only on $S$
and satisfy (\ref{eq:pisketch})? The answer here depends on the question
asked. In this paper, we are trying to prove that the loop-erasure
of random walks on $\delta \ZZ ^{2}$ and $\frac{1}{2}\delta \ZZ ^{2}$
are similar. Therefore we need the graphs $G_{i}$ to be something
that, on certain squares $d+[0,1]^{2}$, $d\in D\subset \ZZ ^{2}$
is similar to $\frac{1}{2}\delta \ZZ ^{2}$ and on others to $\delta \ZZ ^{2}$.
We call such graphs {}``hybrid.'' On a certain intuitive level,
it seems obvious that when we construct this kind of graph the random
walk on it will be similar to Brownian motion, for any defining set
$D$ (or in other words, for any dissection of $\CC $ into squares
of the two types). On a formal level, this requires delicately sawing
together the transition areas (the {}``seams'' in the terminology
of this paper) and lots of technical details. This process is covered
in chapter \ref{chapter_wired_graph}. It starts with the definition
of a hybrid graph and the first step is showing the existence of a
harmonic potential (lemma \ref{lemma_a_wired}). Regrettably, this
particular step requires some computer use, which is described in
the appendix. With the harmonic potential defined, chapter \ref{chapter_wired_graph}
becomes a run-of-the-mill usage of comparisons of continuous and discrete
harmonic functions, and culminates in lemma \ref{lemma_hybrid_square}.
(\ref{eq:pisketch}) is a direct consequence of it, see lemma \ref{lem_lerw_hit_local}.

\subsection{Reading recommendations}

Chapter \ref{chap_generalities} contains various known or unsurprising
facts about random walks and loop-erased random walks. Experts might
want to skip or skim this part. Chapter \ref{chapter_wired_graph},
as explained above, develops the concept of a hybrid graph, a kind
of interpolation between two different graphs, in particular between
two grids of different step length, and shows that the random walk
is not very different from the regular random walk. It is \textbf{highly}
technical and can be skimmed by all. Read carefully the definition
of a hybrid graph, and then the formulation of all lemmas but skip
their proofs. This will not have a significant impact on your ability
to understand later parts. The most interesting part is chapter \ref{chap_core},
with the core being the main lemma, and, to a lesser extent, lemmas
\ref{lem_lerw_hit_local} and \ref{lem_clerw_local}. I recommend
to read it all, linearly, and take a breather after the main lemma.
Starting from section \ref{sec:cont}, the proof gets {}``lighter''
as there is no more need for the machinery of hybrid graphs. All notations
are simpler and techniques are classical. In this part of the proof
(section \ref{sec:cont} and chapter \ref{chap_limit}) the only notable
proof element is lemma \ref{lem_cont_E}. We wrap the proof up in
chapter \ref{chap_limit} which is a two-pages exercise in standard
limit techniques that gives the classical formulation in terms of
the weak limit. It features, though, an interesting example where
loop-erased random walk does \textbf{not} converge (page \pageref{page:exam_punct})
and the exact statement of the theorem (page \pageref{the_theorem}).

\section{\label{chap_generalities}Generalities}

\subsection{Notations}

A weighted graph is a couple $G=(V,W)$ with $V$ the set of vertices
and $W\, :\, V\times V\rightarrow [0,\infty [$, $W(v,w)=W(w,v)$.
Unflinchingly we shall confuse $G$ with $V$, using set notations
such as $v\in G$.

A path in $G$ is a sequence $\gamma =\{\gamma _{i}\}$, $\gamma _{i}\in G$
with $W(\gamma _{i},\gamma _{i+1})\neq 0$. A path is simple if $i\neq j$
implies $\gamma _{i}\neq \gamma _{j}$. The segment of a simple path
$\gamma $ between two points $\gamma _{i}$ and $\gamma _{j}$ is
the subpath $\{\gamma _{i},...,\gamma _{j}\}$ (or the reverse, if
$j<i$). A subset $A\subset G$ is graph-connected if there is a path
in $A$ from every $v\in A$ to every $w\in A$.

For a finite path $\gamma =\{\gamma _{i}\}$ in a graph $G$ we define
its loop erasure, $\LE (\gamma )$, which is a simple path in $G$,
by the consecutive removal of loops from $\gamma $. Formally, \begin{eqnarray*}
\LE (\gamma )_{1} & := & \gamma _{1}\\
\LE (\gamma )_{i+1} & := & \gamma _{j_{i}+1}\quad j_{i}:=\max \{j\, :\, \gamma _{j}=\LE (\gamma )_{i}\}\quad .
\end{eqnarray*}
It will be convenient to consider $\LE (\gamma )$ as a set of vertices
and edges so that we can consider the reversal of $\LE (\gamma )$
as identical to $\LE (\gamma )$, and so that we can write $\LE (\gamma )\cup ...$ 

A random walk on a weighted graph $G$ is a process $R$ that moves
at the $n$th step from $R(n)$ to $R(n+1)$ with the probability
\begin{equation}
\frac{W(R(n),R(n+1))}{\sum _{v}W(R(n),v)}\quad .\label{equ:p_from_W}\end{equation}
If $A\subset B\subset G$ and $v\in G$ then we denote by \[
q(v,A,B,G)\]
the probability of a random walk on $G$ starting from $v$ to hit
$B$ in $A$. A {}``hit'' is only considered for $t\geq 1$ so that
$v\in B$ does not imply a degenerate distribution. If $b\in B$ we
shall write $q(v,b,B,G)$ as a short hand for $q(v,\{b\},B,G)$. 

The Laplacian on a weighted graph $G$ is an operator on functions
$f:G\rightarrow \RR $ (or to any linear space over $\RR $), \[
(\Delta _{G}f)(v)=\sum _{w}W(v,w)(f(w)-f(v))\quad .\]
Clearly, if $T$ is a stopping time for $R$ such that $R(0),...,R(T-1)\not \in B\subset G$,
and $f$ is harmonic (i.e. $\Delta f\equiv 0$) on $G\setminus B$
then\[
\EE f(R(T))=f(R(0))\quad .\]

If $G\subset \CC $ is a graph and $\dreg \subset \CC $, we define\begin{eqnarray*}
\partial _{G}\dreg  & := & \left\{ v\in \dreg \cap G\, :\, \exists w\in G,\, W(v,w)\neq 0\; \wedge \; \left]v,w\right[\not \subset \dreg \right\} \cup \\
 &  & \left\{ w\in G\setminus \dreg \, :\, \exists v\in G,\, W(v,w)\neq 0\; \wedge \; \left]v,w\right[\not \subset \CC \setminus \dreg \right\} \\
\dreg ^{\circ } & := & (G\cap \dreg )\setminus \partial _{G}\dreg 
\end{eqnarray*}
where $\left]v,w\right[$ is the open segment between $v$ and $w$.
We will hardly use the regular definitions of $\partial \dreg $ and
$\dreg ^{\circ }$ so there is little room for confusion. If $v\in \CC $
we define a {}``random walk on $G$ starting from $v$'' as a random
walk on $G$ starting from the point of $G$ closest to $v$. If more
than one exist, choose the top-left point. This also applies to the
notation $q(v,A,B,G)$.

When we say about a set $\dreg \subset \CC $ that it is a polygon
we mean that its boundary is a collection of linear segments of positive
length, but not necessarily that it is simply connected. Punctures
(i.e.~holes of a single point), however, are not allowed.

For a compact metric space $X$, we denote by $\Haus (X)$ the space
of closed subsets of $X$ with the Hausdorff metric, \[
d_{\Haus }(A_{0},A_{1})=\max _{i=0,1}\sup _{a\in A_{i}}d(a,A_{1-i})\]
where as usual $d(b,A)=\inf _{a\in A}d(b,a)$. $\Haus (X)$ is also
a compact metric space. By $\Meas (X)$ we denote the space of measures
on $X$ with the topology of weak convergence. 

$\NN $ denotes the natural integers ($\geq 1$). $\ZZ $ are all
the integers. $\DD $ will denote the disc $|z|<1$ and $\TT $ is
the circle $\partial \DD $. When we write e.g.~$z_{0}+R\DD $ we
mean the usual set addition and multiplication, so it evaluates to
the set $\{z\, :\, |z-z_{0}|<R\}$. The only exception to this rule
is that when $E\subset \RR $ then the notation $E^{2}$ will be used
as a short hand for $E+iE\subset \CC $. In particular, $\ZZ ^{2}$
will be considered as a subset of the complex plane $\CC $ and also
as a graph where \[
W(z,z')=\left\{ \begin{array}{ll}
 1 & |z-z'|=1\\
 0 & \textrm{otherwise}\end{array}
\right.\quad .\]

The notation $\one _{A}$ for a set $A$ stands for the function which
is one on $A$ and zero outside $A$. The support of a function $f$,
denoted by $\supp f$, is the set where $f(x)\neq 0$. The notations
$\wedge $, $\vee $ and $\neg $ are used (somewhat informally) as
shorts for {}``and'', {}``or'' and {}``not'' respectively. The
notation $P\sim Q$ means that the variables $P$ and $Q$ are identically
distributed.

By $C$ and $c$ we denote constants, which could change from formula
to formula (or even inside the same formula).  $C$ will usually pertain
to constants {}``large enough'' and $c$ to constants {}``small
enough''. Occasionally we shall number them for clarity. The notation
$x\approx y$ will be a shorthand for $cy\leq x\leq Cy$.

\subsection{Auxiliary results}

\begin{lem}
\label{lem_wilson}(Wilson's algorithm) The uniform random spanning
tree of a graph $G$ can be constructed using the following inductive
process: in the first step, the partially constructed tree will be
one arbitrary vertex $v\in G$. On the $n$th step $(n>1)$, pick
$w_{n}$ not in the partially constructed tree and add to the latter
a loop-erased random walk on $G$ starting from $w_{n}$ and stopped
when first hitting the partially constructed tree. Continue until
the tree spans all of $G$.
\end{lem}
We do not care what the {}``uniform random spanning tree of $G$''
is (though it is what you would guess). Only that it does not depend
on the algorithm for picking the $v$ and the $w_{n}$'s. This lemma
allows to get all kinds of symmetries for loop-erased random walks,
particularly that the loop-erased random walk from $v$ to $w$ is
distributed identically to the loop-erased random walk from $w$ to
$v$ (though that particular fact was known before).%
\footnote{See lemma \ref{lem_cont_a} for a different use of Wilson's algorithm.%
} The proof can be found in \cite{W96}.

\begin{lem}
\label{lem_condLE_sym}Let $b_{0},b_{1}\in B\subset G$. Let $R_{i}$
be a random walk starting at $b_{i}$, stopped at $B$ and conditioned
to hit $b_{1-i}$. Then \[
\LE (R_{0})\sim \LE (R_{1})\]

\end{lem}
\begin{proof}
Let $R_{i}'$ be a random walk starting at $b_{i}$, stopped at $B$
and conditioned to hit $\{b_{0},b_{1}\}$. Clearly \[
R_{i}\sim R_{i}'\, |\, R_{i}'\textrm{ hits }b_{1-i}\quad .\]
Now let $h$ be the solution of Dirichlet's problem on $G$, with
the initial conditions \[
h(b_{0})=h(b_{1})=1\quad h(B\setminus \{b_{0},b_{1}\})=0\]
 and let $G'$ be a weighted graph with $V'=V$ and $W'(v,w)=h(v)h(w)$.
Let $R_{i}''$ be an (unconditioned) random walk on $G'$ starting
at $b_{i}$ and stopped at $\{b_{0},b_{1}\}$. It is easy to see that
\[
R_{i}'\sim R_{i}''\]
 so \[
R_{i}\sim R_{i}''\, |\, R_{i}''\textrm{ hits }b_{1-i}\quad .\]
 Finally, denoting by $R_{i}'''$ a random walk in $G'$ from $b_{i}$
to $b_{1-i}$ we clearly get \[
\LE (R_{i}''')\sim \LE (R_{i}''\, |\, R_{i}''\textrm{ hits }b_{1-i})\]
and for $\LE (R_{i}''')$ we can use Wilson's algorithm to get \[
\LE (R_{0}''')\sim \LE (R_{1}''')\quad .\qedhere \]

\end{proof}
\begin{lem}
\label{lemma_a_on_Z} \newC{a2ndord}There exists a function $a$
on $\ZZ ^{2}$ such that \begin{eqnarray}
\Delta a & = & \delta _{0}\label{Del_a_del}\\
a(z) & \geq  & a(0)\nonumber \\
a(z) & = & \frac{1}{2\pi }\log |z|+R(z),\quad |R(z)|\leq \frac{C_{\ref {a2ndord}}}{|z|^{2}}\label{asimp_a}
\end{eqnarray}

\end{lem}
A nice proof with a weaker estimate can be found in \cite[section 12.3]{S76}.
The value of $a(0)$ is calculated in \cite[chapter 15]{S76} (note
that Spitzer's $a$ is $4(a+a(0))$ with respect to mine) and is $-\frac{\log 8+2\gamma }{4\pi }$.
A proof that is missing only the actual calculation of $C_{\ref {a2ndord}}$
can be found in \cite{S49} (warning: 60 pages in German). Finally,
see \cite{KS} for a high-order expansion of $a$ and an exact calculation
of $C_{\ref {a2ndord}}=0.017205...$ This function is called the (two
dimensional) discrete harmonic potential.

\begin{lem}
\label{lem_escape_corner2}\newc{cescape} Let $s:=N-1+i(N-d)$ ($1\leq d\leq N$)
and let $I\subset \partial [-1,1]^{2}$ be a connected set with $\diam I>c_{\ref {cescape}}$
and $d(I,(s+1)/N)>c_{\ref {cescape}}$. Then \[
q(s,NI\cap \ZZ ^{2},\partial [-N,N]^{2},\ZZ ^{2})\approx dN^{-2}\]
provided that $N$ is sufficiently large. The constants implicit in
$\approx $ and the minimal $N$ depend on $c_{\ref {cescape}}$. 
\end{lem}
\begin{sublem}Denote $p_{I}(d,N):=q(s,NI\cap \ZZ ^{2},\partial [-N,N]^{2},\ZZ ^{2})$.
Then \[
p_{I}(d,N)\approx dN^{-2}+O\left(N^{-4}\log N\sum _{m=1}^{2N}mp'(m)\right)\]
 where $p'(m):=\max _{d\leq m}p_{H}(d,m)$ and \[
H:=[-1-i,1-i]\cup [-1-i,-1+i]\]
 and the constants implicit in the $\approx $ and in the $O(\cdot )$
above depend on $c_{\ref {cescape}}$. \end{sublem} \begin{proof}[Subproof]Denote
$S=\left]-N,N\right[^{2}$ and denote the value inside the $O(\cdot )$
by $E$. Choose $J\subset \partial [-1,1]^{2}$ to be a connected
set satisfying $I\subset J$ and \[
d(J,(s+1)/N)>\third d(I,(s+1)/N)\qquad d(\partial [-1,1]^{2}\setminus J,I)>\third d(I,(s+1)/N).\]
 Our aim is to prove \[
p_{I}(d,N)\leq CdN^{-2}+CE,\quad p_{J}(d,N)\geq cdN^{-2}-CE\]
which is enough, since we can then exchange the roles of $I$ and
$J$ to get a lower estimate for $p_{I}$.

Let $\varphi $ be the Riemann mapping of $[-1,1]^{2}$ on $\DD $,
$\varphi (0)=0$, $\varphi '(0)>0$. The reflection principle through
the boundary (twice around the corners) for $\varphi $ gives us that
$\varphi $ is analytic near every point of the boundary and in particular
$\varphi _{0}''''\leq C$. From this and from the fact that $\varphi $
preserves the angle near non-corners and doubles the angle near the
corners we get that $\varphi '(b)=0$ only if $b$ is a corner, and
at the corners $\varphi '(b)=0$ and $\varphi ''(b)\neq 0$ --- these
can be summed up as \begin{equation}
\varphi '(b)\approx d(b,K)\label{deriv_phi_is_dist}\end{equation}
 where $K$ is the set of corners, $\{-1,1\}+\{-i,i\}$.

Let $f$ be a real 5 times differentiable function on $\TT $ with
$f(z)=0$ for $\arg z\in \varphi \left(\partial [-1,1]^{2}\setminus J\right)$,
$f(z)=1$ for $\arg z\in \varphi \left(I\right)$ and $0\leq f\leq 1$
and with $f^{(k)}\leq C$, $k=0,\dotsc ,5$ where $C$ depends only
on $c_{\ref {cescape}}$. Extend $f$ to a harmonic function on $\DD $
(we will call the extended function $f$ as well), let $\tilde{f}$
be the complex conjugate of $f$ with $\tilde{f}(0)=0$ and let $F=f+i\tilde{f}$.
It is easy to see that $F^{(k)}\leq C$ for $k=0,1,2,3,4$. We define
a function $\bar{g}$ on $G:=\ZZ ^{2}\cap \bar{S}$ by\[
\bar{g}(v):=f(\varphi (v/N))\quad .\]
Expanding $F(\varphi (v/N))$ to a power series around $v/N$ and
using the fact that \begin{equation}
(z+1)^{k}+(z-1)^{k}+(z+i)^{k}+(z-i)^{k}-4z^{k}=0\quad k=0,1,2,3\label{sym_Z2}\end{equation}
we get (here we used the boundedness of $\varphi ^{(k)}$ and $F^{(k)}$)
\begin{equation}
|\Delta _{G}\bar{g}(v)|\leq CN^{-4}\quad .\label{res_sym_Z2}\end{equation}
We {}``fix'' $\bar{g}$ on $S^{\circ }$ as follows: \[
g(z):=\bar{g}(z)-\sum _{x\in S^{\circ }}\Delta \bar{g}(x)\cdot l_{x}(z)\]
where $l_{x}(z):=a(z-x)-r_{x}(z)$, $a$ is the harmonic potential
from lemma \ref{lemma_a_on_Z}, and $r_{x}$ is the solution of Dirichlet's
problem on $G$ with the conditions \[
r_{x}(z)=a(z-x)\quad z\in \partial S.\]
It is clear from these that $g(z)$ is harmonic on $S^{\circ }$ and
on $\partial S$ we have \[
g(z)=\bar{g}(z)\]
 Next we wish to estimate $g(s)-\bar{g}(s)$. Let $R^{s}$ be a random
walk starting from $s$ and stopped when hitting $\partial S\cup \{x\}$
with some $x\not \in \partial S$ (let $t$ be the stopping time).
(\ref{asimp_a}) gives that $a(z-x)\leq C\log N$ and the boundedness
principle gives the same for $r_{x}(z)$. Because $l_{x}(z)$ is harmonic
on $S^{\circ }\setminus \{x\}$ we get \begin{equation}
l_{x}(s)=\EE (l_{x}(R^{s}(t)))\leq C\log N\cdot \PP (R^{s}(t)=x)\quad .\label{est_axz}\end{equation}
Let $m:=\max |\real (s-x)|,|\imag (s-x)|$. If $m\geq d$ then clearly
$\PP (R^{s}(t)=x)\leq p_{H}(d,m)$. If $m<d$ then a similar argument
gives $\PP (R^{s}(t)=x)\leq 2p_{H}(m,m)$. So in both cases we have
$l_{x}(s)\leq Cp'(m)\log N$ and hence \begin{equation}
|g(s)-\bar{g}(s)|\leq CN^{-4}\log N\sum _{m=1}^{2N}mp'(m)=CE\label{est_g_g_bar}\end{equation}
$\bar{g}(s)$ is easy to estimate (using (\ref{deriv_phi_is_dist}))
because we have \[
d(\varphi (s/N),\TT )\approx N^{-1}\varphi '((s+1)/N)\approx dN^{-2}\]
and an estimate of $f$ using the Poisson kernel and the fact that
$d(\varphi ((s+1)/N),\linebreak [0]\varphi (J))>c$, gives \begin{equation}
|\bar{g}(s)|=|f(\varphi (s/N))|\approx d(\varphi (s/N),\TT )\approx dN^{-2}\quad .\label{est_g_bar}\end{equation}
 Translating the estimates on $g$ to an estimate on the probability
$p(d,N)$ is done by again examining the random walk $R^{s}$ starting
from $s$ but this time stopped on $\partial S$ (let $t$ be the
stopping time). Now $g$ is harmonic on $S^{\circ }$, $g$ is one
on $NI$ and on $\partial S\setminus NI$ we have, $0\leq g\leq 1$.
All these give \begin{eqnarray*}
g(s) & = & \EE (g(R^{s}(t)))\\
 & = & p_{I}(d,N)\EE (g(R^{s}(t))\, |\, R^{s}(t)\in I)+(1-p_{I}(d,N))\EE (g(R^{s}(t))\, |\, R^{s}(t)\not \in I)\\
 & \geq  & p_{I}(d,N)
\end{eqnarray*}
 and similarly $g(s)\leq p_{J}(d,N)$. With (\ref{est_g_g_bar}) and
(\ref{est_g_bar}) the sublemma is proved. \end{proof} \begin{proof}[Proof of lemma \ref{lem_escape_corner2}]
First we estimate $p_{H}$ from above: we use the sublemma, plugging
the estimate $p'(N)\leq 1$ in the right hand side to get \begin{equation}
p_{H}(d,N)\leq CdN^{-2}+CN^{-2}\log N\quad ,\label{eq:pE}\end{equation}
and in particular $p'(N)\leq CN^{-1}$. We now use the sublemma again,
plugging this estimate into the right hand side, and we are done.
\end{proof}

\section{\label{chapter_wired_graph}The hybrid graph}

The proof of the theorem (see page \pageref{main_lemma}) requires
some kind of interpolation between the grids $\frac{1}{N}\ZZ ^{2}$
and $\frac{1}{N'}\ZZ ^{2}$. Before describing the variant I am using,
I wish to make an unusually vague comment. There seems to be some
tradeoff between symmetry and analyticity, in the sense that there
exist models for which it is \textbf{much} easier to prove that the
hybrid process is a good approximation of a Brownian motion, but the
symmetries necessary are not obvious. Being the analyst that I am,
I chose a model for which the proof of lemma \ref{lemma_a_wired}
below is long and technical, but all the symmetries are ready-made
for me. Someone more inclined toward combinatorics might have produced
a nicer proof. 

\begin{defn*}
For a set $D\subset \ZZ ^{2}$, and an integer $N$, we define the
\textbf{hybrid graph} $G(D,N)$, which is a weighted graph, as follows:
The set of vertices $V$ is a union of the following sets: 
\begin{enumerate}
\item The set\[
V_{0}:=\left\{ Nz\in \ZZ ^{2}\, :\, \left\lfloor z\right\rfloor \not \in D\right\} \]
where $\left\lfloor x+iy\right\rfloor :=\left\lfloor x\right\rfloor +i\left\lfloor y\right\rfloor $
i.e.~the vector composed of the two integer values of $x$ and $y$.
\item The set $V_{1}$ which is defined by \[
V_{1}:=\left\{ 2Nz\in \ZZ ^{2}\, :\, \left\lfloor z\right\rfloor \in D,\, N^{-1}\leq d(z,V_{0})\right\} \quad .\]

\end{enumerate}
As for the edges, if $v_{1}$ and $v_{2}\in V_{n}$ then we put an
edge connecting them $\Leftrightarrow $ $|v_{1}-v_{2}|=2^{-n}/N$
and make its weight $1$. If $v_{1}\in V_{0}$, $v_{2}\in V_{1}$
and $|v_{1}-v_{2}|=N^{-1}$ then we connect them by an edge with weight
$\half $ while if $|v_{1}-v_{2}|=N^{-1}\sqrt{\frac{5}{4}}$ then
we connect them by an edge with weight $\quarter $. We denote this
weight by $W$. See figure \ref{fig:hybrid}, left, \vpageref[below]{fig:hybrid}.

The vertices where $V_{0}$ touches $V_{1}$ are called the \textbf{seams}
and are denoted by $\bar{G}$: \[
\bar{G}:=\bigcup _{n=0,1}\{v\in V_{n}\, :\, \exists w\in G\setminus V_{n}\; \wedge \; W(v,w)\neq 0\}\]
We note that $\#\bar{G}\leq CN\#D$ where $\#X$ is the number of
elements of a set $X$. The seams are relevant because we want to
use an argument similar to (\ref{sym_Z2})-(\ref{res_sym_Z2}) on
our hybrid graph. Thus if $f$ is an analytic function we can write
\begin{equation}
\Delta _{G}f(v)\leq C\frac{\max \{f^{(k)}z\, :\, |z-v|<N^{-1}\}}{N^{k}}\label{sym_hybrid}\end{equation}
 where $k=4$ outside $\bar{G}$. On $\bar{G}$, outside the {}``seam-intersections''
we still have (\ref{sym_Z2}) for $i=0$ and $1$, so we can take
$k=2$ (this is easy to verify). Thus we define the seam-intersections
$\bar{\bar{G}}$ using \[
v\in \bar{\bar{G}}\Leftrightarrow \sum _{w}W(v,w)(v-w)\neq 0\quad .\]
 On $\bar{\bar{G}}$ we can only take $k=1$, but luckily there are
even less of these: $\#\bar{\bar{G}}\leq 8\#D$.\label{page:sizesk}
An example of $\bar{G}$ and $\bar{\bar{G}}$ illustrated is in figure
\ref{fig:hybrid}, right.%
\begin{figure}
[b]\input{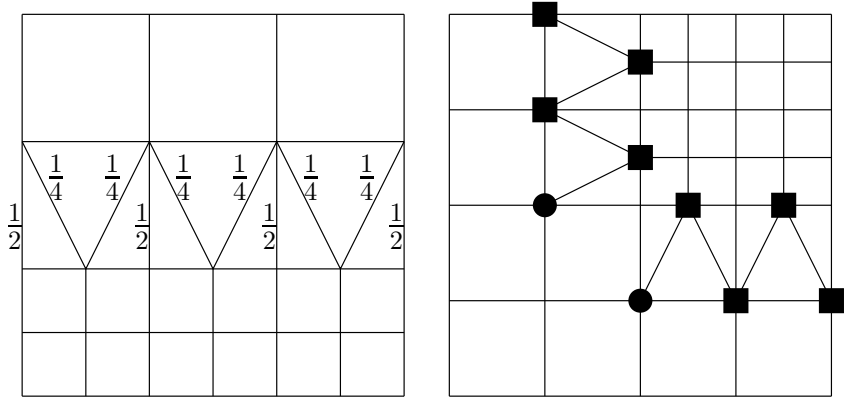}

\caption{\label{fig:hybrid}On the left, a hybrid graph near a seam with the
weights marked (all unmarked edges have weight $1$). On the right,
a hybrid graph near a seam intersection. The set $D$ is $\{0\}$
and $0$ is the middle of the image. Vertices from $\bar{G}$ are
marked with a square and vertices from $\bar{\bar{G}}$ are marked
with a circle.}
\end{figure}
 
\end{defn*}
Eventually, (see page \pageref{random_hybrid}) we sh\-all examine
random walks on a hybrid graph with a random set $D$, so this model
is (locally) a variation on random walk in a random environment. This
might lead the reader to assume that he is in for logarithmic drift
and other cool effects. This is not so --- the model $G$ was constructed
to avoid these effects, and in particular, for all $A$ the drift
is negligible as lemma \ref{lemma_a_wired} and later \ref{lemma_hybrid_brownian}
will demonstrate. 

\begin{lem}
\label{lem:planar}The hybrid graph $G$ is planar. 
\end{lem}
This is easy to verify. \newC{cadmiss} \newC{Cadmisspow} \newc{cadmissres}

\begin{lem}
\label{lemma_a_wired} There exists a $C_{\ref {cadmiss}}$ and a
$C_{\ref {Cadmisspow}}$ such that for any $D\subset [-M,M]^{2}$,
for all $N>C_{\ref {cadmiss}}M^{C_{\ref {Cadmisspow}}}$ the graph
$G=G(D,N)$ satisfies that there exists functions $a_{v}\, :\, G\rightarrow \RR $
with \begin{eqnarray}
\Delta _{G}a_{v}(w) & = & \delta _{v}(w)\label{laplac_a_delta}\\
a_{v}(w) & = & K_{v}\log (N|v-w|)+O\left(N|v-w|\right)^{-c_{\ref {cadmissres}}}\quad \forall v\neq w\label{avw_hybrid_logarithmic}\\
a_{v}(v) & = & O(1)\label{hybrid_ests}\\
K_{v} & \approx  & 1\quad .\label{Kv_1}
\end{eqnarray}
$c_{\ref {cadmissres}}$ depends on $C_{\ref {Cadmisspow}}$.
\end{lem}
The dependency above between $N_{0}$ (the minimal allowed $N$) and
$\diam D$ is not the best possible. The proof can be refined, to
work for certain infinite $D$'s, though not to general ones --- a
checkerboard, i.e. $D=(2\ZZ +2i\ZZ )\cup (2\ZZ +2i\ZZ +1+i)$ seems
to be a particularly bad example. 

\begin{sublem}\label{sublemma_convolution}
Let $G$ be a metric graph and let $b_{v}$ be functions satisfying
\[
\Delta _{G}b_{v}=\delta _{v}+R_{v}\]
 with \[
\sum _{w}|R_{v}(w)|\leq \beta ,\quad \beta <1\]
 and $\bigcup _{v}\supp R_{v}$ finite. Then there exists an $a_{v}$
satisfying (\ref{laplac_a_delta}). Furthermore,

\begin{enumerate}
\item There exist coefficients $\tau _{v,w}$ such that\begin{equation}
a_{v}=\sum _{w}\tau _{v,w}b_{w}\label{a_from_b}\end{equation}
with\begin{equation}
\sum _{w}|\tau _{v,w}|\leq \frac{1}{1-\beta }\label{sum_tau}\end{equation}

\item If $\beta <\frac{1}{2}$ then\begin{equation}
\sum _{w}\tau _{v,w}>\frac{1-2\beta }{1-\beta }\label{sum_tau_large}\end{equation}

\item If\[
\sum _{w\, :\, |v-w|>r}|R_{v}(w)|\leq Ar^{-\mu }\]
then\begin{equation}
\sum _{w\, :\, |v-w|>r}|\tau _{v,w}|\leq C(\beta )Ar^{-\mu }\quad .\label{cvw_small}\end{equation}

\end{enumerate}
\end{sublem}
Of course, the graph $G$ is weighted, but these weights appear only
in the definition of the Laplacian $\Delta _{G}$. The proof below
is a simple use of convolution on $L^{1}$ spaces. \begin{proof}[Subproof]
Define $a_{v}^{1}=b_{v}$ and inductively \begin{equation}
a_{v}^{n}=a_{v}^{n-1}-\sum _{w}(\Delta a_{v}^{n-1}-\delta _{v})(w)a_{w}^{1}\quad .\label{def_a_v_n}\end{equation}
 Notice that the fact that $\bigcup _{v}\supp R_{v}$ is finite gives
that the sum is finite for all $n$. This allows to write \begin{eqnarray}
\sum _{w}|(\Delta a_{v}^{n}-\delta _{v})(w)| & = & \sum _{w}\left|\sum _{x}(\Delta a_{v}^{n-1}-\delta _{v})(x)R_{x}(w)\right|\nonumber \\
 & \leq  & \sum _{x}|(\Delta a_{v}^{n-1}-\delta _{v})(x)|\sum _{w}|R_{x}(w)|\nonumber \\
 & \leq  & \beta \sum _{x}|(\Delta a_{v}^{n-1}-\delta _{v})(x)|\leq \beta ^{n}\label{delta_a_n_v}
\end{eqnarray}
 which gives that $a_{v}^{n}$ converge and that $a_{v}:=\lim _{n\rightarrow \infty }a_{v}^{n}$
satisfies $\Delta a_{v}(w)=\delta _{v}$. (\ref{a_from_b}), (\ref{sum_tau})
and (\ref{sum_tau_large}) are also clear because defining \begin{equation}
\tau _{v,w}:=\delta _{v}(w)-\sum _{n=1}^{\infty }\Delta a_{v}^{n}(w)-\delta _{v}(w)\label{def_tau}\end{equation}
works. We are left therefore with (\ref{cvw_small}). We wish to get
for every $n$ an estimate of the type \[
\sum _{w\, :\, |v-w|>r}|\Delta a_{v}^{n}(w)|\leq A_{n}r^{-\mu }\quad .\]
 Let $d$ satisfy $\beta <d^{\mu }<1$. In the second line of (\ref{delta_a_n_v})
we divide the sum over $x$ as follows: \[
\sum _{w\, :\, |w-v|>r}|\Delta a_{v}^{n}(w)|\leq \beta \bigg (\sum _{|x-v|\leq dr}+\sum _{|x-v|>dr}\bigg )\]
The first sum can be estimated by\[
\sum _{|x-v|\leq dr}|(\Delta a_{v}^{n-1}-\delta _{v})(x)|\sum _{|w-v|>r}|R_{x}(w)|\leq \beta ^{n-1}A_{1}(r-dr)^{-\mu }\]
while the second sum can be estimated by\[
\sum _{|x-v|>dr}|(\Delta a_{v}^{n-1}(x)|\sum _{|w-v|>r}|R_{x}(w)|\leq A_{n-1}(dr)^{-\mu }\beta \]
so we can write\begin{equation}
A_{n}\leq \frac{\beta ^{n-1}A_{1}}{(1-d)^{\mu }}+A_{n-1}\frac{\beta }{d^{\mu }}\quad .\label{An_recurs}\end{equation}
The choice of $d$ is now clear as it ensures that $A_{n}$ converges
exponentially to zero, and in particular $\sum A_{n}<\infty $. This
finishes (\ref{cvw_small}) and the sublemma.\end{proof}

\begin{defn*}
If $s=s_{1}+is_{2}$ and $s_{j}\in \half \ZZ $ we define a function
$A(s,\cdot )$ on $\ZZ ^{2}$ as follows: 
\begin{enumerate}
\item If $s_{1}$ and $s_{2}$ are integers, we take $A(s,v)=a(s-v)$ with
$a$ the harmonic potential on $\ZZ ^{2}$ defined in lemma \ref{lemma_a_on_Z}
above; 
\item If $s_{1}\not \in \ZZ $ and $s_{2}\in \ZZ $, we define \begin{eqnarray*}
t^{\pm } & := & s\pm \half \\
A(s) & := & \frac{1}{2}(A(t^{-})+A(t^{+}))\quad ;
\end{eqnarray*}

\item If $s_{1}\in \ZZ $ and $s_{2}\not \in \ZZ $ we define $A(s)$ symmetrically;
\item If both $s_{i}\not \in \ZZ $ we define \begin{eqnarray*}
t^{\pm ,\pm } & := & s\pm \half \pm \half i\\
A(s) & := & \frac{1}{4}(A(t^{-,-})+A(t^{-,+})+A(t^{+,-})+A(t^{+,+}))\quad .
\end{eqnarray*}

\end{enumerate}
\end{defn*}
\begin{sublem}\label{sublem:est_A}
For all $s$, $\Delta A(s,\cdot )$ is zero except possibly at the
four integer points nearest to $s$, and has the estimate \begin{eqnarray}
A(s,v) & = & \frac{1}{2\pi }\log |s-v|+R_{s}(v)\label{A_logarithmic}\\
|R_{s}(v)| & \leq  & \frac{C}{|s-v|^{2}}\label{est_A_general}
\end{eqnarray}
\end{sublem}
The proof is a simple verification of the 4 cases above and we shall
omit it. 

For the following two sublemmas it will be convenient to use the (somewhat
non-standard) notation $\ZZ ^{+}=\{0,1,...\}$, $\ZZ ^{-}=\ZZ \setminus \ZZ ^{+}$.

\begin{sublem}\label{sublemma_one_seam}
For $D=\ZZ ^{-}+i\ZZ $ we have a function \[
b_{v}(w)\, :\, G(D,N)\rightarrow \RR \]
 such that for every rectangle $S$, \begin{equation}
\sum _{w\in S}|\Delta b_{v}(w)-\delta _{v}(w)|\leq \frac{C}{Nd(v,S\cap \bar{G})}\label{delta_avw_far}\end{equation}
\end{sublem}
\begin{proof}[Subproof]
Define\begin{equation}
b_{v}(w)=\left\{ \begin{array}{ll}
 A(Nv,Nw)-\frac{1}{2\pi }\log N & \left\lfloor w\right\rfloor \not \in D\\
 A(2Nv,2Nw)-\frac{1}{2\pi }\log 2N & \left\lfloor w\right\rfloor \in D\end{array}
\right.\label{def_avw_simp}\end{equation}
 which makes it clear that $\Delta b_{v}-\delta _{v}$ is different
from zero only on $\bar{G}$. Since in our case $\bar{\bar{G}}=\emptyset $
we can use (\ref{sym_hybrid}) with $k=2$, to get for any analytic
function $f$\begin{equation}
|(\Delta _{G}f)(v)|\leq CN^{-2}\max _{|z-v|\leq N^{-1}}|f''(z)|\label{two_deriv}\end{equation}
 and since $\log |z|$ is the real part of such a function (and using
(\ref{A_logarithmic})) we get \[
|\Delta b_{v}(w)|\leq \frac{CN^{-2}}{|v-w|^{2}}\quad .\]
This obviously gives (\ref{delta_avw_far}).
\end{proof}

\begin{sublem}\label{sublemma_two_seams}
Let $D$ be a union of $0-4$ of the quarter planes $\ZZ ^{\pm }+i\ZZ ^{\pm }$.
Then we have a function $b_{v}(w)\, :\, G(D,N)\rightarrow \RR $ with
\begin{equation}
\sum _{w}|\Delta b_{v}(w)-\delta _{v}(w)|\leq 0.4\quad .\label{two_seams_est}\end{equation}
Further, if $S$ is any rectangle then (\ref{delta_avw_far}) also
holds.
\end{sublem}
\begin{proof}[Subproof]
We define $b$ by (\ref{def_avw_simp}). (\ref{delta_avw_far}) follows
easily from sublemma \ref{sublemma_one_seam} and estimating the sum
on $\bar{\bar{G}}$ by (\ref{sym_hybrid}) with $k=1$ and by $\#\bar{\bar{G}}\leq 2$.
(\ref{two_seams_est}) was done numerically and is summed up in appendix
\ref{sect:hyrbid_numerology} (page \pageref{sect:hyrbid_numerology}). \end{proof} \begin{proof}[Proof of lemma \ref{lemma_a_wired}]
The lemma will follow from sublemma \ref{sublemma_convolution} with
the function $b_{v}(w)$ again defined by (\ref{def_avw_simp}). Defining
$r(w)=|\Delta b_{v}(w)-\delta _{v}|$ we need only estimate $\sum r(w)$,
which is non-zero only on $\bar{G}$. Let now $v\in G$. Let $q$
be the integer point closest to $v$ (if more than one exists, choose
any), and let $S(v)=q+[-0.9,0.9]^{2}$ (we may assume $N>10$). (\ref{delta_avw_far})
outside $S$ gives \begin{equation}
\sum _{w\in \bar{G}\setminus S(v)}r(w)\leq CM^{2}/N\quad .\label{sum_r5p}\end{equation}
and with (\ref{two_seams_est}), \[
\sum |(\Delta b_{v}-\delta _{v})(w)|\leq 0.45+CM^{2}/N\]
 so for $N>CM^{2}$ (i.e.~$C_{\ref {Cadmisspow}}=2$) we can use
sublemma \ref{sublemma_convolution} and get (\ref{laplac_a_delta}).
Further, taking $C_{\ref {Cadmisspow}}>2$ we get from (\ref{sum_r5p})
and (\ref{delta_avw_far}) that \begin{eqnarray*}
s<1 & \Rightarrow  & \sum _{|w-v|>s}r(w)\leq \frac{C}{Ns}+CN^{2C_{\ref {Cadmisspow}}^{-1}-1}\leq C\left(Ns\right)^{-c}\\
s\geq 1 & \Rightarrow  & \sum _{|w-v|>s}r(w)\leq \frac{CM^{2}}{Ns}\leq C(Ns)^{-c}
\end{eqnarray*}
from which (\ref{cvw_small}) gives us the same estimate for the coefficients
$\tau $ in (\ref{a_from_b}). To use that, define \[
K_{v}:=\sum _{w}\tau _{v,w}\quad .\]
First note that the requirement (\ref{Kv_1}) follows from (\ref{sum_tau})
and (\ref{sum_tau_large}). Next, (\ref{a_from_b}) gives \[
a_{v}(x)-\frac{1}{2\pi }\log |v-x|\sum _{w}\tau _{v,w}=\sum _{w}\tau _{v,w}\left(b_{w}(x)-\frac{1}{2\pi }\log |v-x|\right)=\Sigma _{1}+\Sigma _{2}\]
where $\Sigma _{1}$ denotes the sum on $w$ satisfying $N|v-w|\leq (N|v-x|)^{1/2}$
and $\Sigma _{2}$ denotes the reminder. For $\Sigma _{1}$, (remember
(\ref{A_logarithmic}))\begin{eqnarray*}
b_{w}(x) & = & \frac{1}{2\pi }\log |w-x|+O\left(N|w-x|\right)^{-2}\\
 & = & \frac{1}{2\pi }\log |v-x|+O\left(\left(N|v-x|\right)^{-1/2}+\left(N|v-x|\right)^{-2}\right)
\end{eqnarray*}
and since $\sum _{w}\tau _{v,w}\leq C$ we get \begin{equation}
\Sigma _{1}=O\left(N|v-x|\right)^{-1/2}\label{x_near_v}\end{equation}
To estimate $\Sigma _{2}$ we write \[
\Sigma _{2}=\sum _{k=1}^{\infty }\Sigma _{2,k}\]
 where $\Sigma _{2,k}$ is the sum on $w$ satisfying $(N|v-x|)^{k/2}<N|v-w|\leq (N|v-x|)^{(k+1)/2}$.
In this case we can estimate \[
\left|b_{w}(x)-\frac{1}{2\pi }\log |v-x|\right|\leq Ck\log N|v-x|\]
 and with $\tau _{v,w}\leq C(N|v-w|)^{-c}$ we get\[
\Sigma _{2,k}\leq C\left(\left(N|v-x|\right)^{k/2}\right)^{-c}Ck\log N|v-x|\]
which we sum over all $k$ and get\begin{equation}
\Sigma _{2}\leq C\left(N|v-x|\right)^{-c}\quad .\label{x_far_from_v}\end{equation}
 (\ref{x_near_v}) and (\ref{x_far_from_v}) give (\ref{avw_hybrid_logarithmic}).
(\ref{hybrid_ests}) is an immediate consequence of (\ref{laplac_a_delta}),
(\ref{Kv_1}) and (\ref{avw_hybrid_logarithmic}). This finishes lemma
\ref{lemma_a_wired}. \end{proof}

\begin{defn*}
A hybrid graph for which $D\subset [-M,M]^{2}$ and $N>C_{\ref {cadmiss}}M^{C_{\ref {Cadmisspow}}}$
is called \textbf{admissible}. We explicitly reiterate the requirement
$C_{\ref {Cadmisspow}}>2$ (which was also used in the proof of lemma
\ref{lemma_a_wired}).
\end{defn*}

\subsection{Global estimates}

In this section will shall prove some simple estimates of hitting
probabilities of random walks on admissible hybrid graphs where the
probability involved is (approximately) independent of $N$. These
are much easier than, for example, estimates for the hitting probability
of a single point, as in lemma \ref{lemma_hybrid_square} further
on.

\begin{lem}
\label{lemma_hybrid_brownian}Let $f$ be a harmonic function on a
domain $E\subset \CC $, $1<\diam E<\infty $. Let $G=G(D,N)$, be
an admissible hybrid graph. Then there exists a function $f'$ on
$G\cap E$, $G$-harmonic on $E^{\circ }$ with \[
|f-f'|\leq CK(f)(\diam E)^{2}N^{-1}\log N\]
with $K(f)$ the maximum on $E$ of all partial derivatives of $f$
up to and including order 4.
\end{lem}
Here, and in other lemmas formulated similarly, we in effect fix the
multiplicative constant in the requirement on $N$ before everything
else, i.e.~the lemma should read {}``There exists some $C$ such
that for all $f$ ...''.

\begin{proof}
By locally adding to $f$ the complex conjugate $\tilde{f}$ we can
use (\ref{sym_hybrid}) to get\[
|(\Delta _{G}f)(v)|\leq CK(f)N^{-k(v)}\]
where \begin{equation}
k(v):=\left\{ \begin{array}{ll}
 4 & v\not \in \bar{G}\\
 2 & v\in \bar{G}\setminus \bar{\bar{G}}\\
 1 & v\in \bar{\bar{G}}\end{array}
\right.\quad .\label{def_kv}\end{equation}
On the other hand, \begin{equation}
\#\{v\in [-M,M]^{2}\, :\, k(v)=k\}\leq \left\{ \begin{array}{ll}
 CM^{2} & k=1\\
 CNM^{2} & k=2\\
 CN^{2}M^{2} & k=4\end{array}
\right.\label{eq:size_Sk}\end{equation}
(see page \pageref{page:sizesk} for these size estimates). This gives
\[
\sum _{v\in E^{\circ }}|(\Delta _{G}f)(v)|\leq CK(f)(\diam E)^{2}N^{-1}\quad .\]
Defining\[
f'=f-\sum _{v\in E^{\circ }}(\Delta _{G}f)(v)\cdot a_{v}\]
we get the required result with (\ref{avw_hybrid_logarithmic}).
\end{proof}
\begin{rem*}
We shall typically use lemma \ref{lemma_hybrid_brownian} to show
that if we have a random walk in a good domain $\dreg $ (smooth boundary)
starting from a point $v$ not too near the boundaries, then the probability
to hit a sizable portion $I$ of the boundary is $>c$. This is done,
as in lemma \ref{lem_escape_corner2}, by taking the solution $f$
of the (continuous) Dirichlet problem on $\dreg $ with $f=1$ on
$I$ and $0$ on $\dreg \setminus I$ (or a smooth approximation of
that), approximating $f$ with a $G$-harmonic $f'$ and using $f'=\mathbb{E}f'(\textrm{hit point})$.
The following lemma is an example.
\end{rem*}
\begin{lem}
Let $z\in \CC $ and $r>1$ some number. Let $R$ be a random walk
on an admissible hybrid graph $G=G(D,N)$, starting from $v$ where
$\frac{3}{4}r<|v-z|<\frac{4}{3}r$. Let $t$ be the stopping time
when $|R(t)-z|<\frac{1}{2}r$ or $|R(t)-z|>2r$. Then \[
\PP \{R|_{[0,t]}\textrm{ contains a loop around }z\}>c\]

\end{lem}
\begin{proof}
Use the previous lemma a number of times (4 should do) to force $R$
to turn $2\pi $ around $z$ and then cross itself. If $N>Cr^{C_{\ref {Cadmisspow}}}$
we can use lemma \ref{lemma_hybrid_brownian} directly. Otherwise
we have $r>CM$ and then on the annulus $A:=(z+2r\DD )\setminus (z+\frac{1}{2}r\DD )$
we have\[
G\cap A=\frac{1}{N}\ZZ ^{2}\cap A\]
 and the question is equivalent to taking $D'=\emptyset $, $r'=1$
and $N'=Nr$, for which we can again use lemma \ref{lemma_hybrid_brownian}.
\end{proof}
\begin{lem}
\label{lem_hybrid_Kesten}\newc{c:kesten} Let $G=G(D,N)$ be an admissible
hybrid graph, let $K\subset G$ be a connected group of vertices,
let $v\in G$, $d(v,K)>1$, and let $R$ be a random walk starting
from $v$ and stopped on $\partial _{G}(v+r\DD )$, $r>\max (1,\diam K)$.
Then \[
\PP \{R\cap K=\emptyset \}\leq C\left(\frac{\diam K}{d(v,K)}\right)^{-c_{\ref {c:kesten}}}\quad .\]

\end{lem}
\begin{proof}
Let $r_{n}=d(v,K)\cdot 2^{n}$, let $T_{k}$ be the stopping times
and $n_{k}$ the numbers defined inductively by \[
T_{k}=\min \{t:\exists n_{k}\neq n_{k-1}\; \wedge \; R(t)\in \partial (v+r_{n_{k}}\DD )\}\quad ;\]
and let $E_{k}$ be the events that $R$ does a loop around $v$ between
$T_{k-1}$ and $T_{k}$. The process does not stop before we have
at least $L\geq c\log (\diam K/d(v,K))$ $E_{k}$'s for which $n_{k}\neq 0$
and the previous lemma gives a lower bound for the probability of
$E_{k}\, |\, \neg E_{1},\ldots ,\neg E_{k-1}$. The fact that the
graph is planar (lemma \ref{lem:planar}) means that the event $E_{k}$
implies that $R$ will necessarily intersect $K$ between $r_{n_{k}-1}$
and $r_{n_{k}+1}$. Therefore \[
\PP \{R\cap K=\emptyset \}\leq \prod _{k=1}^{L}\PP \{\neg E_{k}\, |\, \neg E_{1},\ldots ,\neg E_{k-1}\}\leq c^{L}\quad .\qedhere \]

\end{proof}
\begin{rem*}
For Brownian motion the constant $c_{\ref {c:kesten}}$ is $\frac{1}{2}$
(this is not too difficult to see --- for example, one can use L\"owner's
differential equation to prove that the minimal probability happens
when $K$ is a straight half line and then calculate the probability
explicitly). $c_{\ref {c:kesten}}=\frac{1}{2}$ also for a simple
random walks --- see e.g.~\cite{K87} where an equivalent result
is proved. I have no reason to assume this is not true in our case
too, but we shall not need it. The simple proof above is taken from
\cite[lemma 2.1]{S00}.
\end{rem*}

\subsection{Local estimates}

The aim of this section is to prove lemma \ref{lemma_hybrid_square}
(see also the simplified representation (\ref{est_pbu_2})) which
describes the hitting probability of a point using the geometry of
the domain combined with the local structure of the graph.\newc{chitlog}

\begin{lem}
\label{lem_hit_log} Let $G=G(D,N)$ be an admissible hybrid graph,
let $r>2/N$, let $v\in G$, let $E=\partial _{G}(v+r\DD )\cup \{v\}$,
and let $R$ be a random walk starting from some $w\in G$,$c_{\ref {chitlog}}r<|v-w|<\frac{1}{2}r$
and stopped on $E$. then the probability $p$ that $R$ hits $E$
at $v$ is\[
\approx \frac{1}{\log rN}\quad .\]
The constants implicit in the $\approx $ depend on $c_{\ref {chitlog}}$.
\end{lem}
\begin{proof}
Let $t$ be the stopping time. Let $a_{v}(w)$ be the harmonic potential
from lemma \ref{lemma_a_wired} with respect to the point $v$. Then\[
a_{v}(w)=\EE a_{v}(R(t))=(1-p)\EE (a_{v}(R(t))\, |\, R(t)\neq v)+pa_{v}(v)\]
 and plugging in (\ref{avw_hybrid_logarithmic}) we get \[
p=\frac{K_{v}\log \frac{|v-w|}{r}+O(Nr)^{-c}}{a_{v}(v)-K_{v}\log rN+O(Nr)^{-c}}\]
and since we can assume $r>\frac{C}{N}$ (the possibility to hit $v$
is always positive) and using (\ref{hybrid_ests}) we get $p\approx \log ^{-1}rN$.
\end{proof}
The following definition binds together a number of conditions that
are not really essential but make calculations and proofs easier,
hence the name.

\begin{defn*}
Let $G=G(D,N)$ be an admissible hybrid graph. We say about a rectangle
$S=[r,s]+i[t,u]$ that it is \textbf{easy} in $G$ if $r,s,t,u\in \frac{1}{N}\ZZ $
and $\str S\cdot \diam S>\frac{1}{N}$ where $\str S$ is defined
as the maximal number satisfying\begin{equation}
\str S\leq \frac{r-s}{t-u}\leq \frac{1}{\str S}\label{S_approx_sq}\end{equation}
and\begin{eqnarray}
d(K(S),\bar{G}\cap S^{\circ }) & \geq  & \str S\cdot \diam S\label{S_easy1}\\
d(\partial S,\bar{\bar{G}}\cap S^{\circ }) & \geq  & \str S\cdot \diam S\label{S_easy2}
\end{eqnarray}
and where $K(S)=\{r,s\}+i\{t,u\}$ is the set of corners of $S$.
\end{defn*}
We shall only be interested in rectangles for which $\str S$ is relatively
large. Think about $\str S\geq 0.1$ if you want to get a good notion
of what this definition is all about. Further, when we say {}``let
$S$ be an easy rectangle...'' we always mean in addition {}``with
$\str S$ bounded below by a universal constant''.

\label{page:noteeasy}We note that it follows from (\ref{S_easy1})
and (\ref{S_easy2}) that a seam can only intersect the boundary of
an easy rectangle perpendicularly (see figure \ref{fig:kappa} \vpageref{fig:kappa}).
As this feature seem to follow from local properties of $G$ in a
manner that looks a little random, we note it here. While it simplifies
notations here and there, this is not a significant feature of this
definition.

\begin{lem}
\label{lem_hit_hybrid}\newc{ceasy}Let $G=G(D,N)$ be admissible
and let $S$ be an easy rectangle in $G$. Then for every $u\in S$,
$d(u,\partial S)>c_{\ref {ceasy}}\diam S$ and $b\in \partial S$
we have\begin{eqnarray}
q(u,b,\partial S,G) & \approx  & \frac{d(b,K(S))}{N\diam ^{2}S}\label{hit}\\
q(b,u,\partial S\cup \{u\},G) & \approx  & \frac{d(b,K(S))}{N\diam ^{2}S\log (N\diam S)}\label{escape}
\end{eqnarray}
provided that $N>C+C(\diam S)^{C_{\ref {Cadmisspow}}}$. The constants
implicit in the $\approx $ in (\ref{hit}) and (\ref{escape}) above
depend on $c_{\ref {ceasy}}$ and on $\str S$.
\end{lem}
\begin{proof}
Denote $m=N\diam S$ and $d=Nd(b,K(S))$ so that the right hand sides
of (\ref{hit}) and (\ref{escape}) become $d/m^{2}$ and $d/(m^{2}\log m)$
respectively. Since (\ref{S_easy1}) implies that the hitting probabilities
of the corners are always $0$, we may assume that $m$ is sufficiently
large (the minimal $m$ will depend on $c_{\ref {ceasy}}$ and $\str S$).
We start with (\ref{escape}). Denote $p(b,u):=q(b,u,\partial S\cup \{u\},G)$.
We divide the $b$'s into two cases: \newc{easycorner} \newc{easycircle} \newC{Cforexpo}

\textbf{case 1}: $d<c_{\ref {easycorner}}m$ (we shall fix $c_{\ref {easycorner}}$
later). Assume for simplicity that $b$ is closest to the lower left
corner, $k$ --- the other 3 corners are identical. In this case we
define $\epsilon :=c_{\ref {easycorner}}\diam S$ and \begin{eqnarray*}
Y_{1} & := & k+[0,2\epsilon ]^{2}\\
X_{1} & := & \partial _{G}\left((k+2\epsilon +i[\epsilon ,2\epsilon ])\cup (k+[\epsilon ,2\epsilon ]+2i\epsilon )\right)
\end{eqnarray*}
For $c_{\ref {easycorner}}<\frac{1}{2}\str S$ we have that $G\cap Y_{1}$
is simply a regular grid and we can use lemma \ref{lem_escape_corner2}
to get \[
q(b,X_{1},\partial Y_{1},G)\geq Cdm^{-2}\quad .\]
Next we define \begin{equation}
X_{2}:=\partial _{G}\left(u+\epsilon \DD \right)\quad .\label{last_X}\end{equation}
For $c_{\ref {easycorner}}<\frac{1}{4}c_{\ref {ceasy}}$ we get $d(X_{1},X_{2})>c\diam S$
and then we can use lemma \ref{lemma_hybrid_brownian} to get \[
q(x_{1},X_{2},\partial S\cup X_{2},G)\approx 1\quad \forall x_{1}\in X_{1}.\]
Finally on $X_{2}$ we use lemma \ref{lem_hit_log} to get \[
q(x_{2},u,\partial S\cup \{u\},G)\approx \frac{1}{\log m}\quad \forall x_{2}\in X_{2}\]
 --- we use here \[
x_{2}\in u+4\epsilon \DD \subset S\subset u+C(\diam S)\DD \quad .\]
 This gives $p(b,u)\geq cdm^{-2}\log ^{-1}m$. The estimate of $p(b,u)\leq Cdm^{-2}\brmul \log ^{-1}m$
is identical, but uses a larger $X_{1}$, namely\[
X_{1}:=\partial _{G}Y_{1}\setminus \partial _{G}S\quad .\]
Notice that we are now able to fix $c_{\ref {easycorner}}=\min \left(\frac{1}{3}\str S,\frac{1}{5}c_{\ref {ceasy}}\right)$,
for example.

\textbf{case 2}: $d\geq c_{\ref {easycorner}}m$. This is only slightly
more complicated. Again assume for simplicity that $b$ is in the
lower side of $S$. We start with\begin{eqnarray*}
Y_{1} & := & b+I+i\left[0,\frac{\diam S}{\log ^{2}m}\right]\\
X_{1} & := & \partial _{G}\left(b+I+i\frac{\diam S}{\log ^{2}m}\right)\\
I & := & \left[-C_{\ref {Cforexpo}}\frac{\diam S}{\log m},C_{\ref {Cforexpo}}\frac{\diam S}{\log m}\right]
\end{eqnarray*}
where $C_{\ref {Cforexpo}}$ will be fixed later. To estimate $q(b,X_{1},\partial Y_{1}\cup \partial S,G)$
examine the function $\imag v$. For $m$ sufficiently large the condition
(\ref{S_easy2}) implies $\bar{\bar{G}}\cap Y_{1}=\emptyset $, and
thus $\imag v$ is $G$-harmonic on $Y_{1}$. Let $b'\in S^{\circ }$
be a neighbor of $b$, so $\imag b'-\imag b\approx N^{-1}$, $R$
a random walk starting from $b'$ and stopped on $\partial _{G}Y_{1}$
and let $t$ be the stopping time. Then \[
\imag b'=\EE \imag R(t)\]
 so we get \[
q(b',X_{1},\partial _{G}Y_{1},G)\approx m^{-1}\log ^{2}m\]
 provided we show that the probability to exit $Y_{1}$ on the {}``sides''
(outside $X_{1}\cup \partial S$) is small. But this probability is
clearly (e.g.~by the technique of lemma \ref{lemma_hybrid_brownian})
exponential in the ratio of the length and width of $Y_{1}$ so by
picking $C_{\ref {Cforexpo}}$ sufficiently large we can ignore it.
Summing over all neighbors $b'$ --- usually there is only one but
if $b\in \bar{G}$ there could be two%
\footnote{$b$ cannot have three neighbors in $S^{\circ }$ because the seams
always intersect $\partial S$ perpendicularly --- see the comment
just after the definition of an easy rectangle \vpageref{page:noteeasy}.%
} --- we get \[
q(b,X_{1},\partial _{G}Y_{1}\cup \partial _{G}S,G)\approx m^{-1}\log ^{2}m\quad .\]
Next define \begin{eqnarray*}
Y_{2} & := & b+[-\epsilon ,\epsilon ]+i[0,2\epsilon ]\\
X_{2} & := & \partial _{G}(b+[-\epsilon ,\epsilon ]+2i\epsilon )
\end{eqnarray*}
where again $\epsilon :=c_{\ref {easycircle}}\diam S$ where $c_{\ref {easycircle}}$
is to be defined later, and use lemma \ref{lemma_hybrid_brownian}
to show that, for $m$ sufficiently large, \[
q(x_{1},X_{2},\partial Y_{2},G)\approx \frac{1}{\log ^{2}m}\quad \forall x_{1}\in X_{1}.\]
Finally define $X_{3}$ similarly to (\ref{last_X})\[
X_{3}:=\partial _{G}\left(u+\epsilon \DD \right)\]
and repeat the process of case 1 to get $p(b,u)>cm^{-1}\log ^{-1}m$.
As in case 1, the estimate $p(b,u)<Cm^{-1}\log ^{-1}m$ follows by
merely replacing $X_{2}$ with $\partial Y_{2}\setminus \partial S$.
We see that it is enough to pick $c_{\ref {easycircle}}=\frac{1}{5}c_{\ref {ceasy}}$.

Thus (\ref{escape}) is finished. To get (\ref{hit}) we use the symmetry
of random walk: \begin{equation}
q(u,b,\partial S,G)=\frac{q(b,u,\partial S\cup \{u\},G)}{q(u,\partial S,\partial S\cup \{u\},G)}\cdot \frac{\sum _{v}W(b,v)}{\sum _{v}W(u,v)}\label{def_pab_qa}\end{equation}
where $W$ is the weight function of $G$. we need to estimate $q(u,\partial S,\partial S\cup \{u\},G)$
and we simply sum (\ref{def_pab_qa}) over all $b$ to get \[
q(u,\partial S,\partial S\cup \{u\},G)\approx \sum p(b,a)\approx \frac{1}{\log m}\]
and the lemma is finished.
\end{proof}
\begin{lem}
\label{lem_hit_hybrid2}\newc{chit2}Let $S$ be an easy rectangle
in an admissible hybrid graph $G$. Let $u\in S^{\circ }$ and $b\in \partial S$
with $|u-b|>c_{\ref {chit2}}\diam S$. Then \[
q(u,b,\partial S,G)\approx \frac{d(u,\partial S)d(u,K(S))d(b,K(S))}{N\diam ^{4}S}\]
provided that $N>C+C(\diam S)^{C_{\ref {Cadmisspow}}}$. The constants
implicit in the $\approx $ depend on $c_{\ref {chit2}}$ and on $\str S$.
\end{lem}
The proof is an easy combination of the ideas of the previous proof
(take a square around $u$, a square around $b$, etc) and we shall
omit it.

\begin{lem}
\label{lem:hit_hybrid3} In the previous lemma, without the assumption
$|u-b|>c\diam S$, we get\begin{eqnarray}
q(u,b,\partial S,G) & \approx  & \frac{d(u,\partial S)f(u)f(b)}{N|u-b|^{2}}\label{def_E_lhh3}\\
f(v) & := & \min \left(\frac{d(v,K(S))}{|u-b|},1\right)\quad .\nonumber 
\end{eqnarray}

\end{lem}
\begin{proof}
Let $S_{i}$ be a sequence of easy rectangles, $\str S_{i}\geq c$,
with $u,b\in S_{1}\subset S_{2}\subset \dotsb \subset S_{k}=S$, \[
2\diam S_{i}\leq \diam S_{i+1}\leq C\diam S_{i}\quad ,\]
and where $S_{1}$ satisfies $\diam S_{1}\leq C|u-b|$, $d(u,\partial S_{1})=d(u,\partial S)$,
$d(u,K(S_{1}))\approx f(u)\diam S_{1}$ and $d(b,K(S_{1}))\approx f(b)\diam S_{1}$.
Clearly for some choice of constants such a sequence can always be
found. A little consideration will show that these conditions also
imply $d(b,K(S_{i}))\leq Cf(b)\diam S_{i}$ for all $i$. Define \[
p_{i}:=q(u,b,\partial S_{i},G),\quad q_{i}:=q(u,\partial S_{i}\setminus \partial S,\partial S_{i},G)\quad .\]
Now, lemma \ref{lem_hit_hybrid2} on $S_{1}$ gives $p_{1}\approx E$
and $q_{1}\leq CE_{2}$ where\[
E_{2}:=\frac{d(u,\partial S)f(u)}{|u-b|}\quad E:=E_{2}\frac{f(b)}{N|u-b|}\]
 ($E$ is of course also the right hand side of (\ref{def_E_lhh3})).
For the other $p_{i}$'s we use lemma \ref{lem_hit_hybrid2} on $S_{i+1}$
(and $q_{i}\leq q_{1}\leq CE_{2}$) to get\begin{eqnarray*}
p_{i+1}-p_{i} & = & \sum _{v\in \partial S_{i}\setminus \partial S}q(u,v,\partial S_{i},G)\cdot q(v,b,\partial S_{i+1},G)\\
 & \leq  & C\frac{d(b,K_{i})}{N\diam ^{2}S_{i}}\sum _{v}q(u,v,\partial S_{i},G)\leq C\frac{f(b)}{N\diam S_{i}}E_{2}\quad .
\end{eqnarray*}
which finishes the direction $\leq $ since\[
p_{k}=p_{1}+\sum p_{i+1}-p_{i}\leq CE\left(\sum 2^{-n}\right)=CE\]
 The direction $\geq $ is immediate since $p_{k}\geq p_{1}\geq cE$.
\end{proof}
\begin{lem}
\label{lemma_hybrid_square}\newc{clhs}Let $G=G(D,N)$ be an admissible
hybrid graph. Let $S$ be an easy rectangle in $G$, with $\diam S\approx 1$.
Let $u\in S^{\circ }$ satisfy $d(u,\partial S)>c_{\ref {clhs}}$,
let $b\in \partial S$, and let $p=q(u,b,\partial S,G)$. Then \begin{equation}
p=-\frac{1}{2\pi }\sum _{s\in S^{\circ }}W(b,s)\log |\varphi _{u}(s)|+O\left(pN^{-1/3}\log N\right)\label{hybid_hit_point}\end{equation}
 where $\varphi _{u}$ is the Riemann mapping taking $S$ to $\DD $,
$\varphi _{u}(u)=0$, $\varphi _{u}'(u)>0$, and where $W$ is the
weight function of $G$.

We also assume $N>C(\diam S)^{C_{\ref {Cadmisspow}}}$. The constant
implicit in the $O$ depends on $c_{\ref {clhs}}$, $\str S$ and
on the constants implicit in the condition $\diam S\approx 1$.
\end{lem}
\begin{proof}The proof is based on examining the (unique) solution
$l_{v}$ of the equation \begin{equation} 
\begin{aligned} 
l_{v}(b) & = 0 & b & \in \partial S\\ 
\Delta _{G}l_{v}(z) & = \delta _{v}(z) & z & \in S^{\circ } 
\end{aligned} 
\label{dirichlet} \end{equation} We notice that the maximum principle shows that $l\leq 0$. \begin{sublem}
$|l_{v}(z)|\leq C\log N$ \end{sublem}\begin{proof}[Subproof] We
can write $l_{v}(z)=a_{v}(z)-r_{v}(z)$ where $a_{v}(z)$ comes from
lemma \ref{lemma_a_wired} and $r_{v}(z)$ is the solution of Dirichlet's
problem on $S$ with the conditions $r_{v}(z)=a_{v}(z)$ on $\partial S$.
Lemma \ref{lemma_a_wired} gives that $a_{v}(z)\leq C\log N$ and
the maximum principle gives the same for $r_{v}(z)$. \end{proof}\begin{sublem} \label{sublem_lu_simp}
For $s$ a neighbor of $b\in \partial S$\begin{equation}
|l_{v}(s)|\leq CN^{-1}\log N\min (|v-b|^{-1},d(b,K(S))|v-b|^{-2})\label{lu_simp}\end{equation}
 \end{sublem}\begin{proof}[Subproof] Let $R^{s}$ be a random walk
starting from $s$ and stopped when hitting $\partial S\cup \{v\}$
(let $t$ be the stopping time). Because $l_{v}(z)$ is harmonic on
$S\setminus (\partial S\cup \{v\})$ we get \begin{equation}
-l_{v}(s)=-\EE (l_{v}(R^{s}(t)))\leq C\log N\cdot \PP (R^{s}(t)=v)\label{est_lus_l_bar_us}\end{equation}
 and $\PP (R^{s}(t)=v)\leq Cq(b,v,\partial S\cup \{v\},G)$ can be
estimated using lemma \ref{lem:hit_hybrid3} (use $d(v,\partial S)\leq |v-b|$
and $f(v)\leq 1$) and symmetry (like e.g.~(\ref{def_pab_qa})).
\end{proof} Next some basic facts about $\varphi _{u}$. Denote by
$M$ the middle of $S$. As in the proof of lemma \ref{lem_escape_corner2},
we start with $u=M$, and the reflection principle through the boundary
(twice around the corners) gives us that $\varphi _{M}$ is analytic
near every point of the boundary and in particular $\varphi _{M}''''\leq C$
(here we used the restrictions (\ref{S_approx_sq}) on the geometry
of $S$, and the continuity of the Riemann mapping in the domain%
\footnote{In this case it is easiest to prove this using the Schwarz-Christoffel
formula.%
}). For other $u$ we may take \begin{equation}
\varphi _{u}(x)=F(\varphi _{M}(z)),\quad F(x)=\frac{z-\mu }{1-z\overline{\mu }}\quad .\label{phi_u_is_F_phi}\end{equation}
 where $\mu :=\varphi _{M}(u)$. Explicit differentiation gives $F^{(n)}\leq C$
and hence \begin{equation}
\varphi _{u}(u+z)=A_{u}z+O\left(|z|^{2}\right)\label{asimp_phi_v}\end{equation}
 where $A_{u}:=\varphi '_{u}(u)\approx 1$. Also \begin{equation}
\left|(\log \varphi _{u})^{(n)}(z)\right|\leq \frac{C}{|\varphi ^{n}(z)|}\quad .\label{logphi_hybrid}\end{equation}
\begin{sublem} \label{sublem:lu_deli} Let $s$ be a neighbor of
$b$ and $d=d(b,K(S))N$. Then \begin{equation}
l_{u}(s)=\frac{1}{2\pi }\log |\varphi _{u}(s)|+O\left(dN^{-2}E\right)\label{l_near_boundary}\end{equation}
 where \[
E:=\frac{\log N}{Nr},\qquad r:=\min \left(d(u,\bar{G})+\quarter N^{-1},N^{-2/3}\right)\]
\end{sublem} \begin{proof}[Subproof] Start with the following function
on $G$: let $N'$ be $N$ if $\left\lfloor u\right\rfloor \not \in D$
and $2N$ if $\left\lfloor u\right\rfloor \in D$. Define \begin{equation}
\bar{l}_{u}(v):=\left\{ \begin{array}{ll}
 a(N'(v-u))+\frac{1}{2\pi }\log A_{u}/N' & |v-u|\leq r\\
 \frac{1}{2\pi }\log |\varphi _{u}(v)| & |v-u|>r\end{array}
\right.\label{def_lu_bar}\end{equation}
where $a$ is the harmonic potential on $\ZZ $ --- notice that on
$\{|v-u|\leq d(u,\bar{G})\}$ the hybrid graph is simply $\frac{1}{N'}\ZZ ^{2}$
so the use of $a$ makes sense. In the case $u\in \bar{G}$ (or $r<\frac{1}{N'}$
if you prefer) we simply define $\bar{l}_{u}(u)$ so as to satisfy
$\Delta _{G}\bar{l}_{u}(u)=1$. The uniqueness of $l_{u}$ gives \[
l_{u}(w)=\bar{l}_{u}(w)-\sum _{v\in S^{\circ }\setminus \{u\}}\Delta \bar{l}_{u}(v)\cdot l_{v}(w)\quad .\]
 since the right hand side clearly satisfies (\ref{dirichlet}) (we
shall only use that for $w=s$, though). To estimate $l_{u}-\bar{l}_{u}$
we use sublemma \ref{sublem_lu_simp} for the $l_{u}$'s appearing
on the right hand side and $\Delta \bar{l}_{u}(v)$ is estimated as
follows: 

\begin{enumerate}
\item For $0<|z|<r-N^{-1}$ we have $\Delta \bar{l}_{u}(z+u)=0$.
\item At the transition annulus $|z|=r+O(N^{-1})$, (\ref{asimp_a}) gives
us (for $|z|\leq r$)\begin{equation}
\left|\bar{l}_{u}(z+u)-\frac{1}{2\pi }\log |A_{u}z|\right|\leq \frac{C}{(Nr)^{2}}\label{equ:lu_a_log}\end{equation}
while using (\ref{asimp_phi_v}) together with the fact $A_{u}\geq c$
gives for $r<|z|\leq r+\frac{1}{N}$, \begin{equation}
\left|\bar{l}_{u}(z+u)-\frac{1}{2\pi }\log |A_{u}z|\right|\leq Cr\leq CN^{-2/3}\leq \frac{C}{(Nr)^{2}}\quad .\label{equ:lu_phi_log}\end{equation}

\item Finally, for $|v-u|>r+1$ expand $\log |\varphi |=\real \log \varphi (z)$
to a power series around $z$ and (\ref{sym_hybrid}) will give (using
(\ref{logphi_hybrid}))\[
|\Delta \bar{l}_{u}(v)|\leq \frac{C}{(N|\varphi (v)|)^{k(v)}}\]
($k(v)$ from (\ref{def_kv})).
\end{enumerate}
This division into three cases, combined with the different possibilities
for $k(v)$ and $|v-b|$ is formalized by dividing $G\cap S=\cup _{i=0}^{6}F_{i}$
with the $F_{i}$'s defined as follows: \[
F_{0}:=\left\{ v\, :\, |v-u|<r-N^{-1}\right\} ,\quad F_{1}:=\left\{ v\, :\, \big ||v-u|-r\big |\leq N^{-1}\right\} \quad ;\]
$F_{2}$-$F_{4}$ are subsets of \[
V:=\left\{ v\in S^{\circ }\, :\, r+\frac{1}{N}<|v-u|\; \wedge \; |v-b|>\alpha \right\} ,\quad \alpha :=\half \min \str S\cdot \diam S,c_{\ref {clhs}}\quad ,\]
$F_{2}:=V\setminus \bar{G}$, $F_{3}:=V\cap (\bar{G}\setminus \bar{\bar{G}})$
and $F_{4}:=V\cap \bar{\bar{G}}$; and finally $F_{5}$-$F_{6}$ are
related similarly to $V':=\{|v-b|\leq \alpha \}$ --- there is no
$F_{7}$ because $V'\cap \bar{\bar{G}}=\emptyset $ due to (\ref{S_easy2}).
This gives \begin{eqnarray*}
\left|l_{u}(s)-\log |\varphi _{u}s|\right| & \leq  & \Sigma _{1}+\Sigma _{2}+\Sigma _{3}+\Sigma _{4}+\Sigma _{5}+\Sigma _{6}\\
\Sigma _{i} & := & \sum _{v\in F_{i}}|\Delta \bar{l}_{u}(v)|\cdot |l_{v}(s)|
\end{eqnarray*}
($\Sigma _{0}$ being 0). We now estimate them one by one.

For $\Sigma _{1}$ we use (\ref{equ:lu_a_log}), (\ref{equ:lu_phi_log}),
sublemma \ref{sublem_lu_simp} and $\#F_{1}\leq CNr$ to get \[
\Sigma _{1}\leq \sum _{v\in F_{1}}\frac{C}{(Nr)^{2}}\cdot \frac{Cd\log N}{N^{2}|v-b|^{2}}\leq CdN^{-2}E\quad .\]
For $\Sigma _{2}$ we use the easy fact that $|\varphi _{u}(v)|\geq c|v-u|$
to sum by distance from $u$ and get \begin{eqnarray*}
\Sigma _{2} & \leq  & \sum _{v\in F_{2}}\frac{C}{(N|\varphi (v)|)^{4}}\cdot \frac{Cd\log N}{N^{2}|v-s|^{2}}\leq \frac{Cd\log N}{N^{2}}\sum _{k=Nr}^{N}k\frac{C}{k^{4}}\leq \\
 & \leq  & C\frac{d\log N}{N^{2}}\cdot \frac{1}{(Nr)^{2}}\leq CdN^{-2}E\quad .
\end{eqnarray*}
A similar estimate for $\Sigma _{3}$ and $\Sigma _{4}$ (remember
(\ref{eq:size_Sk})) gives \[
\Sigma _{3},\Sigma _{4}\leq CdN^{-2}E\]
Next we tackle $\Sigma _{5}$. This time we sum by distance from $b$:\begin{eqnarray*}
\Sigma _{5} & \leq  & \frac{C\log N}{N^{4}}\bigg (\sum _{\substack{ v\in F_{5}\brk N|v-b|>d}
}\frac{Cd/N^{2}}{|v-b|^{2}}+\sum _{\substack{ v\in F_{5}\brk N|v-b|\leq d}
}\frac{C/N}{|v-b|}\bigg )\\
 & \leq  & \frac{C\log N}{N^{4}}\left(\sum _{k=d}^{cN}k\frac{Cd/N^{2}}{(k/N)^{2}}+\sum _{k=1}^{d}k\frac{C/N}{k/N}\right)\leq \frac{Cd\log ^{2}N}{N^{4}}\\
 & \leq  & CdN^{-2}E
\end{eqnarray*}
and similarly\[
\Sigma _{6}\leq \frac{C\log ^{2}N}{N^{2}}\]
and since (\ref{S_easy1}) implies that $\Sigma _{6}$ is relevant
only when $d>cN$ this is also $\leq CdN^{-2}E$. Summing the estimates
for $\Sigma _{1},\dotsc ,\Sigma _{6}$ gives us (\ref{l_near_boundary}).
\end{proof} With sublemma \ref{sublem:lu_deli} proved, we are capable
of proving lemma \ref{lemma_hybrid_square} for the case $d(u,\bar{G})>N^{-2/3}$,
and to get some estimate for the other case. We again use the time-symmetry
of random walks in the form (\ref{def_pab_qa}). For the nominator
of (\ref{def_pab_qa}), examine the random walk $R^{s}$ starting
from $s$ and stopped on $\partial S\cup \{u\}$ and $t$ the stopping
time and get \begin{equation}
l_{u}(s)=\EE (l_{u}(R^{s}(t)))=l_{u}(u)\cdot q(s,u,\partial S\cup \{u\},G)\quad .\label{equ:nom}\end{equation}
 so remembering (\ref{equ:p_from_W}),\[
q(b,u,\partial S\cup \{u\},G)=\frac{\sum _{s\in S^{\circ }}W(b,s)l_{u}(s)}{l_{u}(u)\cdot \sum _{s}W(b,s)}\]
For the denominator, we examine a random walk starting from $v$ a
neighbor of $u$ and find in the same manner\[
q(v,\partial S,\partial S\cup \{u\},G)=1-\frac{l_{u}(v)}{l_{u}(u)}\]
and summing over all $v$ we get \begin{equation}
q(u,\partial S,\partial S\cup \{u\},G)=1-\sum _{v}\frac{W(u,v)}{\sum _{w}W(u,w)}\cdot \frac{l_{u}(v)}{l_{u}(u)}=\frac{-(\Delta l_{u})(u)}{l_{u}(u)\cdot \sum _{w}W(u,w)}\label{equ:denom}\end{equation}
so (\ref{equ:nom}), (\ref{equ:denom}) and (\ref{def_pab_qa}) with
(\ref{dirichlet}) give \begin{equation}
q(u,b,\partial S,G)=-\sum _{s}W(b,s)l_{u}(s)\label{equ:p_and_l}\end{equation}
and with (\ref{l_near_boundary}) and $p\approx dN^{-2}$ we get \begin{equation}
p=-\frac{1}{2\pi }\sum _{s\in S^{\circ }}W(b,s)\log |\varphi _{u}(s)|+O\left(pE\right)\quad .\label{est_with_dist_Gbar}\end{equation}
 This proves the lemma for the case $d(u,\bar{G})>N^{-2/3}$.

For the other case, let $u$ satisfy $d(u,\bar{G})\leq N^{-2/3}$.
Let $m\approx N^{-1/2}$ satisfy that the square $S'$ of side length
$m$ around $u$ is easy (clearly such an $m$ can be found). With
this $S'$ we can write \begin{equation}
p=\sum _{v\in \partial S'}q(u,v,\partial S',G)\cdot q(v,b,\partial S,G)\quad ;\label{p_is_prod_qq}\end{equation}
$q(v,b,\partial S,G)$ can be estimated by (\ref{est_with_dist_Gbar})
to give \begin{equation}
q(v,b,\partial S,G)=-\frac{1}{2\pi }\sum _{s\in S^{\circ }}W(b,s)\log |\varphi _{v}(s)|+O\left(pE_{v}\right)\quad ;\label{with_dist_q_to_Gbar}\end{equation}
Thus we have to estimate $\varphi _{u}-\varphi _{v}$. But if $\nu :=\varphi _{u}(v)$
then $|\nu |\leq Cm$ and furthermore \[
\varphi _{v}=F(\varphi _{u}),\quad F(z):=\frac{z-\nu }{1-z\bar{\nu }}\quad .\]
Writing the Taylor expansion of $\log |\varphi _{u}|$ near $b$ and
plugging in the derivatives of $F$ will give \begin{eqnarray*}
\lefteqn{\big ||\log \varphi _{v}(s)|-|\log \varphi _{u}(s)|\big |\leq } &  & \\
 & \leq  & \left|(b-s)\real \frac{\varphi _{u}'}{\varphi _{u}}(b)\right|\cdot O(|\nu |)+\left|\frac{(b-s)^{2}}{2}\real \frac{\varphi _{u}''\varphi _{u}-\varphi _{u}'^{2}}{\varphi _{u}^{2}}(b)\right|\cdot O(|\nu |)+O(N^{-3})\\
 & \leq  & O\left(\frac{|\nu |\cdot |\varphi '|}{N}+\frac{|\nu |}{N^{2}}+N^{-3}\right)
\end{eqnarray*}
so we can replace $\varphi _{v}$ with $\varphi _{u}$ in (\ref{with_dist_q_to_Gbar})
to get \[
q(v,b,\partial S,G)=-\frac{1}{2\pi }\sum _{s\in S^{\circ }}W(b,s)\log |\varphi _{u}(s)|+O\left(pE_{v}+dN^{-5/2}\right)\]
and summing over (\ref{p_is_prod_qq}) we get\begin{equation}
p=-\frac{1}{2\pi }\sum _{s\in S^{\circ }}W(b,s)\log |\varphi _{u}(s)|+O(\Sigma _{1}+\Sigma _{2})\label{eq:p_error_sig1sig2}\end{equation}
with $\Sigma _{1}$ and $\Sigma _{2}$ defined by \begin{eqnarray}
\Sigma _{1} & := & p\sum _{v\in \partial S'}q(u,v,\partial S',G)\cdot E_{v}\nonumber \\
\Sigma _{2} & := & \sum q(u,v,\partial S',G)\cdot dN^{-5/2}=dN^{-5/2}\leq CpN^{-1/2}\label{eq:sig2_phiuphiv_diff}
\end{eqnarray}
To estimate $\Sigma _{1}$ use lemma \ref{lem_hit_hybrid} and get\[
q(u,v,\partial S',G)\leq \frac{C}{m}\]
so summing by distance from $\bar{G}$ we get\begin{eqnarray*}
\Sigma _{1} & \leq  & CpN^{-1/3}\log N+p\sum _{k=1}^{N^{1/3}}\frac{C}{m}\cdot \frac{C\log N}{Nk}\\
 & \leq  & CpN^{-1/3}\log N
\end{eqnarray*}
which with (\ref{eq:p_error_sig1sig2}) and (\ref{eq:sig2_phiuphiv_diff})
proves the lemma. \end{proof}\begin{rmks}

\begin{enumerate}
\item The following weaker form of the lemma will probably look more familiar:\begin{equation}
p=\kappa _{b}\frac{|\varphi _{u}'(b)|}{2\pi N}+O\left(pN^{-1/3}\log N+N^{-2}\right)\label{est_pbu_2}\end{equation}
where the structure constant $\kappa _{b}$ is defined by%
\begin{figure}
[b]\input{kappa.xfig.pstex_t}

\caption{\label{fig:kappa}$\kappa _{b}$ on the boundary of an easy rectangle.}
\end{figure}
 \begin{equation}
\kappa _{b}:=\left\{ \begin{array}{ll}
 \frac{1}{2} & b\in V_{1}\\
 1 & b\in V_{0}\setminus \bar{G}\\
 \frac{9}{8} & b\in V_{0}\cap \bar{G}\end{array}
\right.\quad ,\label{def_kappab}\end{equation}
$V_{0}$ and $V_{1}$ from the definition of a hybrid graph. See figure
\ref{fig:kappa}. To get (\ref{est_pbu_2}) just write a Taylor expansion
of $\log |\varphi |=\real \log \varphi $ near $b$ and a few orientation
arguments will allow to calculate $\arg (b-s)\varphi '/\varphi $.
For example, for $b\in G\setminus \bar{G}$, we get $\arg (b-s)\varphi '/\varphi =0$.
As already remarked, the conditions (\ref{S_easy1}) and (\ref{S_easy2})
imply that $\partial S$ can only intersect a seam perpendicularly
--- otherwise we would need a number of additional special values
for $\kappa _{b}$. The additional error $N^{-2}$ in (\ref{est_pbu_2})
is the second term in the Taylor expansion. This error is of course
meaningful only for $b$ close to the boundaries.
\item The $\log $ factor can be removed. Basically one has to take the
estimates for $l$ given by sublemma \ref{sublem:lu_deli} which,
for $|u-s|>c$ are better than those of sublemma \ref{sublem_lu_simp},
and plug them right back into the estimates of $\Sigma _{1}$-$\Sigma _{4}$.
Also the assumption $\diam S\approx 1$ is unnecessary --- without
it the lemma holds with $m:=N\diam S$ instead of $N$.
\item The final statement of lemma \ref{lemma_hybrid_square} can be translated
back into $l_{u}$ terms using (\ref{equ:p_and_l}) to give an estimate
in (\ref{l_near_boundary}) that does not depend on $d(u,\bar{G})$.
\item The value of $r$ is not the best for any $u$. For example, if $S$
is a square around $u$ then the symmetry of the situation shows that
\[
\varphi (u+z)=-i\varphi (u+iz)=-\varphi (u-z)=i\varphi (u-iz)\]
from which we may conclude \begin{equation}
\varphi (u+z)=A_{u}z+O(|z|^{5})\quad .\label{asymp_phi2}\end{equation}
The improved estimate in (\ref{equ:lu_phi_log}) allows to pick $r=N^{-1/3}$
(assume for simplicity $d(u,\bar{G})>N^{-1/3}$) and to get in (\ref{hybid_hit_point})
the error estimate $O(pN^{-2/3}\log N)$.
\item This is actually a rather nice result even for random walks on $\frac{1}{N}\ZZ ^{2}$
(which is of course a trivial hybrid graph --- take $D=\emptyset $
and get $\kappa _{b}\equiv 1$). For example, it shows that for the
square $[-1,1]^{2}$ the hitting probability from $0$ of $b$ is
$\frac{1}{2\pi N}|\varphi '(b)|+O(N^{-5/3}\log N)$. In comparison,
the probability that a Brownian motion will hit an interval of length
$\frac{1}{N}$ around $b$ is $\frac{1}{2\pi }|\varphi '(b)|+O(N^{-2})$.
In other words, because there are no {}``quantization effects''
we get an error $N^{-2/3}\log N$ better than what we would expect
for, say, a quantized circle.
\item Forgetting for the moment hybrid graphs we reread the proof for the
case of $u$ in the center of a square in $\frac{1}{N}\ZZ ^{2}$.
The role of the symmetries of the grid $\ZZ ^{2}$ seems to suggest
that an equivalent calculation for a random walk on a triangular grid
will give stronger results. Let $p(b)$ be the probability that a
random walk on a triangular grid with step length $\frac{1}{N}$ starting
from 0 will hit a regular hexagon centered at 0 of side length $1$
at the point $b$. We see that we get a better estimate in (\ref{asymp_phi2}),
better estimates in (\ref{asimp_a})%
\footnote{See \cite{KS} for a proof that $a_{T}(z)=\frac{1}{\pi \sqrt{12}}\log |z|+\beta +O(|z|^{-4})$
where $a_{T}$ is the harmonic potential of the triangular grid (notice
that our $a_{T}$ is $\frac{1}{6}$ of the $a$ in \cite{KS}) which
can be normalized to get $\beta =0$. Interestingly, the value $\frac{1}{\pi \sqrt{12}}$
may also be deduced from the proof above by summing (\ref{hybid_hit_point})
(with the factor $\frac{1}{2\pi }$ replaced with the value we are
calculating) over all $b$ and using the fact that $\int |\varphi '|=2\pi $.%
} and $k(v)=6$ for all $v$, so we should be able to get $p=\frac{1}{2\pi N}|\varphi '(b)|+O(N^{-2})$,
no? However, for a hexagon, the inverse Schwarz-Christoffel function
$\varphi $ is no longer analytic around the corners, which gives
an additional error term $O(d^{-1}N^{-5/2}\log N)$ where $d$ is
the distance to the nearest corner. Indeed, at the very corner $\varphi '(z)$
is $0$, while the hitting probability is $\geq cN^{-3/2}$.
\item Trying the formulation $p=-\frac{1}{\pi \sqrt{12}}\sum W(b,s)\log |\varphi (s)|$
in the previous remark improves the error estimates in the hexagon's
edges' middle parts but not near the corners. A careful calculation
will give in this case that the best $r$ is $N^{-0.4}$ and the error
is \[
O(pN^{-1.8}+d^{-5/2}N^{-4})\]
and again, when $d<C/N$ the error becomes $O(N^{-3/2})$, which is
exactly the magnitude of $p$.
\end{enumerate}
\end{rmks}

\section{\label{chap_core}The proof core}

For simplicity of notation, assume throughout this chapter that $M>1$
and $N>1$ (so we don't have to worry about $\log $'s being zero).

\subsection{Localization}

Lemma \ref{lemma_hybrid_square} gave a relatively precise estimate
of the difference between random walk on $\ZZ ^{2}$ and on a hybrid
graph. In this section we mold this general lemma into some corollaries
in the form required for the proof of the theorem. Specifically, we
estimate the amount a random walk changes when the graph is changed
on one square from $\frac{1}{N}\ZZ ^{2}$ to $\frac{1}{2N}\ZZ ^{2}$,
or, in our notation, when one $z\in \ZZ ^{2}$ is moved into or out
of $D$.

\begin{lem}
\label{lem_exer}Assume for $i=1,2$ \[
Q_{i}=\prod _{k=1}^{K}M_{i,k}\]
and for all $k$, $|M_{1,k}-M_{2,k}|\leq \epsilon \min M_{1,k},M_{2,k}$.
Assume $\epsilon K\leq \frac{1}{2}$. Then\[
|Q_{1}-Q_{2}|\leq C\epsilon K\min Q_{1},Q_{2}\quad .\]

\end{lem}
This exercise is left for the reader.

\begin{lem}
\label{lem_lerw_hit_local}Let $D_{1},D_{2}\subset \ZZ ^{2}\cap \left[-M,M\right[^{2}$
with $D_{1}\bigtriangleup D_{2}=\{z\}$ ($\bigtriangleup $ being
the symmetric difference). Let $S:=z+[-1,2]^{2}$. Let $G_{i}=G(D_{i},N)$
be admissible, $N>CM^{C_{\ref {Cadmisspow}}}$. Let $B\subset (G_{1}\setminus S^{\circ })\cap [-M,M]^{2}=(G_{2}\setminus S^{\circ })\cap [-M,M]^{2}$
be a set containing a loop around a point $a\in \CC $ and let $b\in B$
be some point. Then\begin{eqnarray*}
|p_{1}-p_{2}| & \leq  & CN^{-1/3}\log N\log ^{2}M\min \left(p_{1},p_{2}\right)\\
p_{i} & := & q(a,b,B,G_{i})\quad .
\end{eqnarray*}

\end{lem}
\begin{proof}
Let $m_{2}=\frac{1}{N}\left\lfloor \frac{2}{3}N\right\rfloor $, $m_{3}=\frac{1}{N}\left\lfloor \frac{1}{3}N\right\rfloor $
and let $S_{j}=\left[-m_{j},1+m_{j}\right]^{2}$ ($j=2,3$). We note
that both squares are easy for both $G_{i}$. \begin{sublem} \label{sublem_S3_to_S2}
Let $x\in S_{3}$ and let $x_{i}$ be a point of $G_{i}$ closest
to $x$. Let $w\in \partial S_{2}$ and let $q_{i}=q(x_{i},w,\partial S_{2},G_{i})$.
Then\begin{equation}
|q_{1}-q_{2}|\leq CN^{-1/3}\log N\min \left(q_{1},q_{2}\right)\label{del_S3_to_S2}\end{equation}
\end{sublem} \begin{proof}[Subproof] If $x_{1}=x_{2}$ this follows
immediately from lemma \ref{lemma_hybrid_square} because the list
of neighbors $s$ of $w$, $W(w,s)$ and $\varphi _{x}$ do not depend
on $i$. Otherwise we have to estimate $\log |\varphi _{x_{1}}|-\log |\varphi _{x_{2}}|$
using $|x_{1}-x_{2}|\leq 1/N$. For this we use the representation
(coming from expanding the $\log \varphi $ in (\ref{hybid_hit_point})
into a Taylor series)\[
q_{i}=\kappa _{w}\frac{|\varphi '|}{N}+\sum _{s}W(w,s)\real \frac{(s-w)^{2}}{2}\cdot \frac{\varphi ''\varphi -\varphi '^{2}}{\varphi ^{2}}+O(N^{-3})\quad ,\]
where $\kappa _{w}$ is defined in (\ref{def_kappab}) (and is independent
of $i$). $|x_{1}-x_{2}|\leq 1/N$ gives \[
\left|\varphi _{x_{1}}^{(j)}-\varphi _{x_{2}}^{(j)}\right|\leq CN^{-1}\min \left(\varphi _{x_{1}}^{(j)},\varphi _{x_{2}}^{(j)}\right)\qquad j=0,1,2\]
and the sublemma is proved. \end{proof} \begin{sublem} \label{sublem_S2_to_b}
Let $v,w\in \partial S_{2}$, then the probabilities $q^{v}=q(v,b,\partial S_{3}\cup B,G_{i})$
(which obviously don't depend on $i$) satisfy\begin{equation}
q^{v}\approx q^{w}\label{qv_approx_qw}\end{equation}
\end{sublem} \begin{proof}[Subproof] If $d(v,w)<\frac{1}{12}$ this
is easy because one can find some $2d(v,w)<m<\frac{1}{3}$ such that
$S_{4}:=v+\left[-m,m\right]^{2}$ is easy and write\[
q^{x}=\sum _{y\in \partial S_{4}}q(x,y,\partial S_{4},G_{i})\cdot q(y,b,\partial S_{3}\cup B_{i},G_{i})\]
and since $q(x,y,\partial S_{4},G_{i})\approx d(y,K(S_{4}))/N$ for
both $x=v$ and $w$, we are done. If $d(v,w)\geq \frac{1}{12}$ just
write a chain $v_{i}\in \partial S_{2}$, $v_{0}=v$, $v_{n}=w$ and
$|v_{i}-v_{i+1}|<\frac{1}{12}$ and use the first case inductively
($n\leq 17$, for example). \end{proof} \begin{sublem} \label{sublem_hit_b_close}
Let $v\in \partial S_{2}$. then the probabilities $q_{i}=q(v,b,B,G_{i})$
satisfy\begin{equation}
|q_{1}-q_{2}|\leq CN^{-1/3}\log N\log ^{2}M\cdot \min \left(q_{1},q_{2}\right)\label{q1_approx_q2}\end{equation}
\end{sublem} \begin{proof}[Subproof] Let $R_{i}^{v}$ be random walks
on $G_{i}$ starting from $v$. Define stopping times $t_{i}(0)=0$
and\begin{eqnarray*}
t_{i}(2k+1) & := & \min \{t>t_{i}(2k)\, :\, R_{i}(t)\in \partial S_{3}\cup B\}\\
t_{i}(2k) & := & \min \{t>t_{i}(2k-1)\, :\, R_{i}(t)\in \partial S_{2}\}
\end{eqnarray*}
and let $k_{i}$ be the first $k$ such that $R_{i}^{v}(t_{i}(k))\in B_{i}$
($k_{i}$ is always odd, of course). Let $p_{i,k}(w)$ be the probability
that $k_{i}>k$ and $R_{i}^{v}(t_{i}(k))=w$. Then we can write $p_{i,k}(w)$
as a sum \[
p_{i,k}(w)=\sum _{\vec{w}}\PP \left(R_{i}^{v}(t_{i}(l))=w_{l},\; l=0,\dotsc ,k\right)\]
where the sum is over all vectors $\vec{w}=\left\{ w_{l}\right\} _{l=0}^{k}$
with $w_{0}=v$, $w_{k}=w$, $w_{2l}\in \partial S_{2}$ and $w_{2l+1}\in \partial S_{3}$.
Lemma \ref{lem_exer} and sublemma \ref{sublem_S3_to_S2} will now
give for every $k<K:=\lfloor cN^{1/3}\log ^{-1}N\rfloor $, \begin{equation}
|p_{1,k}(w)-p_{2,k}(w)|<CkN^{-1/3}\log N\min \left(p_{1,k}(w),p_{2,k}(w)\right)\quad .\label{pikw_similar}\end{equation}
From these we derive an equivalent estimate for\[
q_{i,k}:=\PP \left(\left\{ k_{i}=k\right\} \cap \left\{ R_{i}^{v}(t_{i}(k))=b\right\} \right)\]
since\begin{equation}
q_{i,k+1}=\sum _{w\in \partial S_{2}}p_{i,k}(w)\cdot q(w,b,\partial S_{3}\cup B,G_{i})\label{qik_sum_pikw}\end{equation}
 and we can sum over $w$. On the other hand, summation over sublemma
\ref{sublem_S2_to_b} shows that the probability to hit $b$ is approximately
independent from $k$, i.e. \begin{equation}
q_{i}\approx \frac{q_{i,2k+1}}{\PP (k_{i}>2k)}\label{hit_b_indep_k}\end{equation}
and a simple exit probability estimate shows that for $k$ odd \begin{equation}
\PP \left(k_{i}=k\, |\, k_{i}\geq k\right)\geq \frac{c}{\log M}\quad .\label{k_exponent}\end{equation}
To get (\ref{k_exponent}) notice that the fact that $B$ contains
a loop around $a$ allows to bound the hitting probability of $B$
by the hitting probability of $\partial [-M,M]^{2}$ which can be
estimated by $\frac{c}{\log M}$, e.g.~using lemma \ref{lemma_hybrid_brownian}.
This gives \begin{eqnarray}
|q_{1}-q_{2}| & \leq  & \sum _{k=1}^{\infty }|q_{1,k}-q_{2,k}|\leq \nonumber \\
 & \leq  & \sum _{k=1}^{K}CkN^{-1/3}\log N\left(1-\frac{c}{\log M}\right)^{k}\min q_{1},q_{2}+\nonumber \\
 &  & +\sum _{k=K+1}^{\infty }C\left(1-\frac{c}{\log M}\right)^{k}\min q_{1},q_{2}\label{q1_approx_q2_calc}
\end{eqnarray}
and we are done. \end{proof} The lemma is now easy: if $a\in S_{3}$
we write \[
q(a,b,B,G_{i})=\sum _{v\in \partial S_{2}}q(a,v,\partial S_{2},G_{i})\cdot q(v,b,B,G_{i})\]
and use sublemma \ref{sublem_S3_to_S2} to show that $q_{i}:=q(a,v,\partial S_{2},G_{i})$
satisfy \[
|q_{1}-q_{2}|\leq CN^{-1/3}\log N\min q_{1}.q_{2}\]
 and get the result using sublemma \ref{sublem_hit_b_close}. If $a\not \in S_{3}$
we similarly write \[
q(a,b,B,G_{i})=q(a,b,B\cup \partial S_{3},G_{i})+\sum _{v\in \partial S_{3}}q(a,v,B\cup \partial S_{3},G_{i})\cdot q(v,b,B,G_{i})\]
and since both $q(a,b,B\cup \partial S_{3},G_{i})$ and $q(a,v,B\cup \partial S_{3},G_{i})$
are independent of $i$, the lemma follows from sublemma \ref{sublem_hit_b_close}
like the previous case.
\end{proof}
\begin{lem}
\label{lem_clerw_local} With the notations of lemma \ref{lem_lerw_hit_local}
and $a\not \in S$, let $\check{R}_{i}$ be a random walk on $G_{i}$
starting from $a$ and conditioned to hit $B$ at $b$. Let $\check{\gamma }_{i}$
be the segment of $\LE (\check{R}_{i})$ until $S$, or all of $\LE (\check{R}_{i})$
if $\LE (\check{R}_{i})\cap S=\emptyset $. Then\[
\sum _{\gamma }|\PP (\check{\gamma }_{1}=\gamma )-\PP (\check{\gamma }_{2}=\gamma )|\leq CN^{-1/3}\log ^{2}M\log ^{3}N\]
where the sum is taken on all the simple paths $\gamma $ in $G_{i}$
from $a$ to $S\cup b$.
\end{lem}
We note that (due to $a\not \in S$) we can assume $a\in G_{1}$ and
get $a\in G_{2}$ too. Note also and that a path in $G_{1}$ from
$a$ to $S\cup b$ is also a path in $G_{2}$ from $a$ to $S\cup b$.

\begin{proof}
We keep all notations from the proof of lemma \ref{lem_lerw_hit_local}.
Denote by $R_{i}^{a}$ a random walk on $G_{i}$ starting from $a$
and stopped on $B$ and let $k_{i}^{a}$ and $t_{i}^{a}$ be the equivalents
(for $R_{i}^{a}$) of $k_{i}$ and $t_{i}$ from sublemma \ref{sublem_hit_b_close}.
Combining (\ref{hit_b_indep_k}) and (\ref{k_exponent}) and summing
over $v\in \partial S_{2}$ we get \begin{eqnarray*}
\PP \left(\left\{ k_{i}^{a}>k\right\} \cap \left\{ R_{i}^{a}\, \mathrm{hits}\, b\right\} \right) & \leq  & C\left(1-\frac{c}{\log M}\right)^{k}\PP \left(\left\{ k_{i}^{a}>1\right\} \cap \left\{ R_{i}^{a}\, \mathrm{hits}\, b\right\} \right)\\
 & \leq  & C\left(1-\frac{c}{\log M}\right)^{k}\PP \left(\left\{ R_{i}^{a}\, \mathrm{hits}\, b\right\} \right)\quad .
\end{eqnarray*}
This means that by taking $K=\lfloor C\log M\log N\rfloor $ we can
write\begin{equation}
\PP (k_{i}^{a}>K\, |\, R_{i}^{a}\, \mathrm{hits}\, b)<CN^{-2}\quad .\label{eq:not_after_K}\end{equation}
Next, denote by $\gamma _{i}$ the unconditioned version of $\check{\gamma }_{i}$
i.e. the segment of $\LE (R_{i}^{a})$ until $S$. Since $\gamma _{i}$
obviously depends only on the portions of $R_{i}^{a}$ outside $S_{2}$
we can, as in sublemma \ref{sublem_hit_b_close}, sum over all vectors
$\vec{w}=\left\{ w_{j}\right\} _{j=1}^{k-1}$, $w_{2l}\in \partial S_{2}$,
$w_{2l+1}\in \partial S_{3}$ and get, using (\ref{del_S3_to_S2})
and lemma \ref{lem_exer},\begin{eqnarray*}
|p_{1,k}-p_{2,k}| & \leq  & CkN^{-1/3}\log N\min \left(p_{1,k},p_{2,k}\right)\\
p_{i,k} & := & \PP \left(\left\{ \gamma _{i}=\gamma \right\} \cap \left\{ k_{i}^{a}=k\right\} \cap \left\{ R_{i}^{a}(t_{i}^{a}(k))=b_{i}\right\} \right)
\end{eqnarray*}
 and summing over $k$ from $1$ to $K$ we get \begin{eqnarray}
|p_{1}-p_{2}| & < & CK^{2}N^{-1/3}\log N\min p_{1},p_{2}+E_{1}+E_{2}\label{p1_approx_p2}\\
p_{i} & := & \PP \left(\left\{ \gamma _{i}=\gamma \right\} \cap \left\{ R_{i}^{a}\, \mathrm{hits}\, b_{i}\right\} \right)\nonumber \\
E_{i} & := & \PP \left(\left\{ \gamma _{i}=\gamma \right\} \cap \left\{ R_{i}^{a}\, \mathrm{hits}\, b_{i}\right\} \cap \left\{ k_{i}^{a}\geq K\right\} \right)\nonumber 
\end{eqnarray}
which finishes the proof since $\check{p}_{i}:=\PP (\check{\gamma }_{i}=\gamma )$
is given by \[
\check{p}_{i}=\frac{p_{i}}{\PP (R_{i}^{a}\, \mathrm{hits}\, b)}\]
and the difference between the nominators can be estimated by lemma
\ref{lem_lerw_hit_local}. The factors $E_{i}$ are dealt with by
(\ref{eq:not_after_K}) which gives $\sum _{\gamma }\frac{E_{i}}{\PP (R_{i}^{a}\, \mathrm{hits}\, b)}\leq CN^{-2}$.
\end{proof}
\begin{defn*}
Let $\gamma $ be a path in a metric graph $G$, let $z\in \ZZ ^{2}$
and let $r>\varepsilon >0$. An $(r,\epsilon ,z)$-\textbf{quasi loop}
of $\gamma $ are two points $v,w\in \gamma $ , $|v-z|\leq \epsilon $,
$|w-z|\leq \epsilon $ such that the section of $\gamma $ between
$v$ and $w$ has a diameter $>r$. We denote $\gamma \in \QL (r,\epsilon ,z)$.
\end{defn*}
Note that this is slightly different than the $(r,\epsilon ,z)$-quasi
loops of \cite{S00}.

\begin{lem}
\label{lem_no_quasi_loops} \newc{cquasi} There exists a constant
$c_{\ref {cquasi}}$ such that for every $G=G(D,N)$ admissible hybrid
graph; every $a\in G$; every $B\subset G\cap [-M,M]^{2}$ containing
a loop around $a$; and every $\epsilon \geq 1$ and $r>M^{1-c_{\ref {cquasi}}}$
we have \begin{equation}
\EE (\#\{z:\LE (R^{a})\in \QL (r,\epsilon ,z)\})<CM^{-c}\epsilon ^{C}\log N\label{eq:no_ql}\end{equation}
where $R^{a}$ is a random walk on $G$ starting from $a$ and stopped
on $B$.
\end{lem}
This lemma follows from lemma \ref{lem_hybrid_Kesten} like lemma
3.4 in \cite{S00} follows from lemma 2.1 ibid. However, the differences
(especially those resulting from the fact that $B$ is not necessarily
the boundary of a simply connected domain) seem to merit a reproduction
of Schramm's proof.

\begin{proof}
We may assume $64\epsilon <r<\sqrt{2}M$. For $4\epsilon <s<\frac{1}{16}r$
and $z\in \CC $ we define the sets $S_{3}:=z+(s+\epsilon )\DD $
and $S_{2}:=z+2s\DD $. As in lemma \ref{lem_clerw_local}, we define
stopping times $t^{a}(k)$ by \begin{eqnarray}
t^{a}(2k+1) & := & \min \{t>t^{a}(2k):R^{a}(t)\in \partial S_{3}\cup B\}\nonumber \\
t^{a}(2k) & := & \min \{t>t^{a}(2k-1):R^{a}(t)\in \partial S_{2}\cup B\}\label{eq:def_tak}
\end{eqnarray}
and $k^{a}$ to be the first $k$ such that $R^{a}(t^{a}(k))\in B$.
Next we define the variable \[
X:=X(z):=\#\{x\in \ZZ ^{2}\cap (z+s\DD ):\LE (R^{a}[0,t^{a}(k^{a})])\in \QL (r,\epsilon ,x)\}\quad .\]
With these notations we can write\begin{equation}
\EE X\leq M^{2}\PP (k^{a}>K)+\sum _{k=1}^{K}\EE (X\cdot \one _{\{k^{a}=k\}})\quad .\label{eq:X_sum_Xk}\end{equation}
and as above (\ref{k_exponent}) holds for $k$ odd, so for $K>C\log ^{2}M$
for some $C$ sufficiently large we get \[
\PP (k^{a}>K)\leq \left(1-\frac{c}{\log M}\right)^{K/2}\leq M^{-C}\]
for any desired constant in the exponent of $M$. In other words,
with this $K$ we can ignore the first summand in (\ref{eq:X_sum_Xk}).
Let us therefore define the number of quasi-loops up to the $k$th
time\[
X_{k}:=\#\{x\in \ZZ ^{2}\cap (z+s\DD ):\LE (R^{a}[0,t^{a}(k)])\in \QL (r,\epsilon ,x)\}\quad .\]
The process of loop-erasing between $t^{a}(2k)$ and $t^{a}(2k+1)$
can only destroy $(r,\epsilon ,x)$-quasi loops for $x\in z+s\DD $
so we get \[
X\cdot \one _{\{k^{a}=2k+1\}}\leq X_{2k}\cdot \one _{\{k^{a}=2k+1\}}\leq X_{2k}\cdot \one _{\{k^{a}>2k\}}\]
 With this in mind we define \[
\Delta _{k}:=(X_{2k+2}-X_{2k})\cdot \one _{\{k^{a}>2k+2\}}\quad .\]
 \begin{sublem} \label{sublem:onegam}Let $v\in \partial S_{3}$
be some vertex; and let $\gamma $ be a simple path on $G$ starting
and ending on $\partial S_{2}$; let $R^{v}$ be a random walk starting
from $v$ and stopped at $\partial S_{2}\cup B$; and define \begin{eqnarray*}
\delta '(v,\gamma ) & := & \#\{x\in \ZZ ^{2}\cap (z+s\DD ):d(x,\gamma )\leq \epsilon \; \wedge \; d(x,R^{v})\leq \epsilon \}\cdot \one _{\{R^{v}\cap \gamma =\emptyset \}}\quad ,\\
\delta (v,\gamma ) & := & \delta '(v,\gamma )\cdot \one _{Y}
\end{eqnarray*}
where $Y$ is the event that $R^{v}$ stops on $\partial S_{2}$.
Then \[
\EE \delta '(v,\gamma )\leq C\epsilon ^{2}\qquad \EE \delta (v,\gamma )\leq Cs^{-c}\epsilon ^{C}.\]
\end{sublem} \begin{proof}[Subproof] \newC{C:lem18t} Define the event
$E(t)$ by\[
E(t):=\{\exists x\in \ZZ ^{2}\cap (z+s\DD ):d(x,\gamma )\leq \epsilon \; \wedge \; d(x,R^{v}(t))\leq \epsilon \}\]
and use $E$ to define stopping times $s_{1},...$ by \[
s_{1}:=\min \{t:E(t)\; \vee \; R^{v}(t)\in \partial S_{2}\cup B\}\]
 and for $i>1$ \[
s_{i}:=\min \{t>s_{i-1}:(E(t)\; \wedge \; |R^{v}(t)-R^{v}(s_{i-1})|>C_{\ref {C:lem18t}}\epsilon )\; \vee \; R^{v}(t)\in \partial S_{2}\cup B\}\]
where $C_{\ref {C:lem18t}}$ will be defined promptly. Let $l$ be
the number of $s_{i}$'s defined before the process is stopped i.e.~$s_{l+1}\in \partial S_{2}\cup B$.
Lemma \ref{lem_hybrid_Kesten} gives for $i<l$ that the probability
not to intersect $\gamma $ between $s_{i}$ and $s_{i+1}$ is $\leq c<1$
assuming $C_{\ref {C:lem18t}}$ is large enough --- fix $C_{\ref {C:lem18t}}$
to satisfy that. This gives $\PP (l>A)\leq c^{A}$ for any value of
$A$. Further, for the time period between $s_{A}$ and $s_{l+1}$
lemma \ref{lem_hybrid_Kesten} gives that \[
\PP (Y\cdot \one _{\{R^{v}\cap \gamma =\emptyset \}}\, |\, l\geq A)\leq q(R^{v}(s_{A}),\partial S_{2},\partial S_{2}\cup \gamma ,G)\leq C(s/\epsilon )^{-c}\quad .\]
Finally, it is clear that $\delta (v,\gamma )>CA\epsilon ^{2}$ implies
$l>A$ therefore\begin{align}
\EE \delta (v,\gamma ) & \leq \sum _{A=1}^{\infty }\EE (\delta (v,\gamma )\cdot \one _{\{l=A\}})\leq \nonumber \\
 & \leq C\epsilon ^{2}\sum _{A=1}^{\infty }A\PP (Y\cap \{R^{v}\cap \gamma =\emptyset \}\cap \{l>A-1\})\leq \nonumber \\
 & \leq C\epsilon ^{2}\sum _{A=1}^{\infty }Ac^{A-1}\cdot C(s/\epsilon )^{-c}\leq Cs^{-c}\epsilon ^{C}\quad .\label{eq:Xvg_Y}
\end{align}
The inequality for $\EE \delta '(v,\gamma )$ follows similarly. \end{proof}
Returning to the estimate of $\Delta _{k}$ we notice that in order
to have $X_{2k+2}>X_{2k}$ we need to have for some $x\in \ZZ ^{2}\cap (z+s\DD )$
that $\LE (R^{a}[0,t^{a}(2k+2)])$ contains an $(r,\epsilon ,x)$-quasi
loop and $\LE (R^{a}[0,t^{a}(2k)])$ doesn't contain one. This requires
at least that
\begin{enumerate}
\item $x$ is $\epsilon $-near at least one segment $\gamma $ of $\LE (R^{a}([0,t^{a}(2k)]))\cap S_{2}$.
\item $R_{k}:=R^{a}([t^{a}(2k+1),t^{a}(2k+2)])$ gets $\epsilon $-near
$x$ and then fails to intersect at least one of the segments $\gamma $
from 1. 
\end{enumerate}
In other words, the number of such $x$'s can be estimated by $\delta (R^{a}(t^{a}(2k+1)),\gamma )$.
With this in mind we denote by $\Gamma (t,u)$ the collection of connected
components $\gamma $ of $\LE (R^{a}([0,t]))\cap (z+u\DD )$ satisfying
$\gamma \cap \partial S_{3}\neq \emptyset $ and $R^{a}(t)\not \in \gamma $
and get \[
\EE \Delta _{k}\leq \sum _{\gamma \in \Gamma _{k}}\max _{v\in \partial S_{3}}\EE \delta (v,\gamma )\leq C(\#\Gamma _{k})s^{-c}\epsilon ^{C}\]
where $\Gamma _{k}:=\Gamma (t^{a}(2k+1),2s)$. It easy to see that
$\#\Gamma _{k}\leq k$, and summing up to $k$ we get\[
\EE X_{2k}\cdot \one _{\{k^{a}>2k\}}\leq Cs^{-c}k^{2}\epsilon ^{C}\quad .\]
Another summation, up to $K$, will give us \begin{eqnarray}
\EE X & \leq  & M^{-C}+Cs^{-c}\epsilon ^{C}\log ^{6}M+\EE X'\label{eq:EX_not_boundary}\\
X' & := & \Delta _{(k^{a}-2)/2}'\cdot \one _{\{k^{a}\text {\, even}\}}\nonumber \\
\Delta _{k}' & := & (X_{2k+2}-X_{2k})\cdot \one _{\{k^{a}>2k+1\}}\quad .\nonumber 
\end{eqnarray}
Thus we are left with the estimate of $\EE X'$, which is the behavior
near the boundary --- if $B\cap S_{2}=\emptyset $ then of course
$k^{a}$ is always odd and we get $X'\equiv 0$. It is at this point
that we utilize the difference between $r$ and $s$. Further, it
will be easier to use entry probabilities rather than exit probabilities.
Thus the first step will be a time-reversed lemma \ref{lem_hybrid_Kesten}.
\begin{sublem} Let $S_{0}:=z+\tsh {\frac{1}{2}}r\DD $ and $S_{1}:=z+\frac{1}{4}r\DD $.
Let $\gamma $ be a simple path between $\partial S_{0}$ and $\partial S_{3}$.
Then \begin{equation}
q(v,\partial S_{3},\partial S_{3}\cup \partial S_{0}\cup \gamma ,G)\leq C(r/s)^{-c}\quad \forall v\in \partial S_{1}.\label{eq:qwv_large}\end{equation}
\end{sublem}
\begin{proof}
[Subproof]Let $\gamma ':=\gamma \cap \left(\left(z+\frac{1}{8}r\DD \right)\setminus S_{2}\right)$
(remember the condition $s<\frac{1}{16}r$) and let $w$ be an arbitrary
point with $\frac{1}{4}s<|w-z|<\frac{3}{4}s$. Lemma \ref{lem_hybrid_Kesten}
shows that \[
q(w,\partial \left(z+\tsh {\frac{1}{8}}r\DD \right),\partial \left(z+\tsh {\frac{1}{8}}r\DD \right)\cup \gamma '\cup \partial (z+\tsh {\frac{1}{4}}s\DD ),G)\leq C(r/s)^{-c}\quad .\]
It is easy to get from that, using lemmas \ref{lemma_hybrid_brownian}
and \ref{lem_hit_log} as in the proof of lemma \ref{lem_hit_hybrid}
that for $w\in \partial (z+\frac{1}{2}s\DD )$ \begin{equation}
q(w,v,\partial S_{0}\cup \gamma '\cup \{w,v\}\cup \partial (z+\tsh {\frac{1}{4}}s\DD ),G)\leq C\frac{(r/s)^{-c}}{\log sN\log rN}\leq C\frac{(r/s)^{-c}}{\log ^{2}N}\label{eq:qvwS0S3gam}\end{equation}
and the symmetry of random walk (in the form (\ref{def_pab_qa}))
gives the same estimate for $q(v,w,\partial S_{0}\cup \gamma '\cup \{w,v\}\cup \partial (z+\frac{1}{4}s\DD ),G)$.
Reversing the argument used to get (\ref{eq:qvwS0S3gam}) we get \[
q(v,\partial S_{3},\partial S_{3}\cup \partial S_{0}\cup \gamma ',G)\leq C(r/s)^{-c}\quad \forall v\in \partial S_{1}.\]
and then of course it holds for $\gamma $ as well.
\end{proof}
\begin{sublem} \label{sublem:isolated}For every $z$, $s$ and $r$,\begin{equation}
\EE X'\leq C\epsilon ^{2}\log ^{4}M\log N\left(\frac{r}{s}\right)^{-c}q(a,B\cap (z+4s\DD ),B,G)\label{eq:sublem_isolated}\end{equation}
\end{sublem} 
\begin{proof}
[Subproof] Let $k>0$ be some integer, and, with the same $S_{0}$
and $S_{1}$ as above, define times $s_{1}$ and $s_{2}$ by\begin{eqnarray*}
s_{1} & := & \max \{t^{a}(1)<t<t^{a}(2k+1):R^{a}(t)\in \partial S_{0}\}\\
s_{2} & := & \min \{t>s_{1}:R^{a}(t)\in \partial S_{1}\}
\end{eqnarray*}
--- if $R^{a}([t^{a}(1),t^{a}(2k+1)])\cap \partial S_{0}=\emptyset $
define both to be $t^{a}(2k+1)$. We notice that $R^{a}$ from $s_{1}$
(and therefore from $s_{2}$ as well) to $t^{a}(2k+1)$ is a random
walk conditioned not to hit $\partial S_{0}$. Lemma \ref{lemma_hybrid_brownian}
gives \[
q(v,\partial S_{3},\partial S_{0}\cup \partial S_{3},G)\geq \frac{c}{\log r/s}\quad \forall v\in \partial S_{1}\]
so (\ref{eq:qwv_large}) gives for any path $\gamma $ from $\partial S_{0}$
to $\partial S_{3}$ \[
\PP (R^{a}([s_{2},t^{a}(2k+1)])\cap \gamma =\emptyset \, |\, R^{a}([0,s_{2}]))\leq C(r/s)^{-c}\log r/s\leq C(r/s)^{-c}\]
and in particular this is true for $\gamma \in \Gamma (s_{2},\frac{1}{2}r)$
so \[
\EE \left.\left(\#\Gamma \left(t^{a}(2k+1),\tsh {\frac{1}{2}}r\right)\; \right|\, k^{a}>0\right)\leq C(r/s)^{-c}\cdot \max \#\Gamma (s_{2},\tsh {\frac{1}{2}}r)\leq Ck(r/s)^{-c}\quad .\]
Denoting $\Gamma _{k}':=\Gamma (t^{a}(2k+1),\frac{1}{2}r)$ we get\[
\EE \Delta _{k}'\leq \EE \sum _{\gamma \in \Gamma _{k}'}\max _{v\in \partial S_{3}}\delta _{k}'(v,\gamma )\leq C\epsilon ^{2}k\left(\frac{r}{s}\right)^{-c}\PP (k^{a}>0)\]
and summing over $k$ we get\[
\EE X'\leq \sum _{k=1}^{\infty }\EE \Delta _{k}'\leq M^{-C}+C\epsilon ^{2}\log ^{4}M\left(\frac{r}{s}\right)^{-c}\PP (k^{a}>0)\quad .\]
The only thing left is to notice that (\ref{eq:sublem_isolated})
is obvious when $B\cap S_{2}=\emptyset $ so we can assume it has
at least one point. This implies that the probability to hit $B\cap (z+4s\DD )$
when starting from an arbitrary point in $\partial S_{3}$ before
exiting from $z+4s\DD $ is $\geq \frac{c}{\log sN}\geq \frac{c}{\log N}$.
Therefore $\PP (k^{a}>0)\leq C\log N\cdot q(a,B\cap (z+4s\DD ),B,G)$
and the sublemma is proved.
\end{proof}
Lemma \ref{lem_no_quasi_loops} now follows by summing over $z$.
Let $z_{1},\dotsc ,z_{l}$ be the points of $s\ZZ ^{2}\cap [-M,M]^{2}$,
so $l\leq C(M/s)^{2}$. Then $[-M,M]^{2}\subset \bigcup z_{i}+s\DD $
and therefore (\ref{eq:EX_not_boundary}) gives \newc{cquasihalf}
\begin{align*} 
\EE (\#\{z:\LE (R^{a})&\in \QL (r,\epsilon ,z)\}) \leq \sum _{i=1}^{l}\EE X(z)\leq \\ 
& \leq C\left( \frac{M}{s}\right) ^{2}s^{-c}\epsilon ^{C}\log ^{6}M+\sum _{i=1}^{l}\EE X'(z)
\end{align*}
which makes it clear that for some $c_{\ref {cquasihalf}}$, taking
$s=M^{1-c_{\ref {cquasihalf}}}$ would make the first summand $\leq C\epsilon ^{C}M^{-c}$.
For the sum on $\EE X'$ we use (\ref{eq:sublem_isolated}) to get\begin{eqnarray*}
\sum _{i=1}^{l}\EE X'(z) & \leq  & C\epsilon ^{2}\log ^{4}M\log N\left(\frac{r}{s}\right)^{-c}\sum _{i=1}^{l}q(a,B\cap (z_{i}+4s\DD ),B,G)\\
 & \leq  & C\epsilon ^{2}\log ^{4}M\log N\left(\frac{r}{s}\right)^{-c}\cdot C
\end{eqnarray*}
so by taking $c_{\ref {cquasi}}=\frac{1}{2}c_{\ref {cquasihalf}}$
we get that the second summand is $\leq C\epsilon ^{C}M^{-c}\log N$
and the lemma is finished.
\end{proof}
\begin{rem*}
As in \cite{S00}, this result (in the case $D=\emptyset $, i.e.~a
regular random walk) implies that for every open bounded set $\dreg $,
every subsequence limit of the random walks on $G:=\delta \ZZ ^{2}$
starting from some $a\in \dreg $ and stopped on $\partial _{G}\dreg $
is supported on the set of simple paths (this follows from lemma \ref{lem_no_quasi_loops}
exactly like theorem 1.1 in \cite{S00} follows from lemma 3.4 ibid.).
Thus we get a strengthening of the second statement of the above mentioned
theorem 1.1 --- it is now true for any open set $\dreg $, without
the restriction that the diameter of every component of $\partial \dreg $
is positive. The example on page \pageref{page:exam_punct} shows
that in this setting the formulation using subsequence limits is necessary
as the limit does not necessarily exists.
\end{rem*}

\subsection{Symmetry}

\begin{mainlem} \label{main_lemma} \newc{cmain} \newC{Cmain} Let
$\dreg \subset [-1,1]$ be an open polygon; let $\delta >0$ and for
$i=1,2$ let $G_{i}=2^{1-i}\delta \ZZ ^{2}$. Let $a\in \dreg $ and
let $a_{i}$ be the point of $G_{i}$ closest to $a$. Let $R_{i}$
be a random walk in $G_{i}$ from $a_{i}$ to $\partial _{G_{i}}\dreg $.
Let $\ereg \subset \dreg $ be some open set, $a\in \ereg $. Then
\begin{eqnarray}
\PP (\LE (R_{2})\subset \ereg +C_{\ref {Cmain}}\delta ^{c_{\ref {cmain}}}\DD ) & > & \PP (\LE (R_{1})\subset \ereg )-C\delta ^{c}\label{R2_geq_R1}\\
\PP (\LE (R_{1})\subset \ereg +C_{\ref {Cmain}}\delta ^{c_{\ref {cmain}}}\DD ) & > & \PP (\LE (R_{2})\subset \ereg )-C\delta ^{c}\label{R1_geq_R2}
\end{eqnarray}
for $\delta <\delta _{0}(\dreg )$. The constants $C_{\ref {Cmain}}$
and $c_{\ref {cmain}}$ are independent of $\dreg $, $\ereg $, $\delta $
and $a$. \end{mainlem} 
\renewcommand{\thesublem}{M.\arabic{sublem}}
\setcounter{sublem}{0}

\begin{proof}
We begin with the proof of (\ref{R2_geq_R1}). Let $M$ be some value
--- we shall fix the best value for $M$ later. Let $N=\frac{1}{M\delta }$.
One of the conditions on $M$ will be that $N\in \NN $. Define \[
\mu :=N^{-1/3}\log ^{3}N\log ^{2}M\]
and \[
p_{1}:=\PP (\LE (R_{1})\subset \ereg )\quad .\]

\textbf{step 1}: Define subsets $Y\subset X\subset \ZZ ^{2}$ as follows:\[
X:=\{x:x+[0,1]^{2}\cap M\dreg \neq \emptyset \}\]
\[
Y:=\left\{ x\in X\, :\, d(x,M(\partial \ereg \setminus \partial \dreg ))>M^{1-c_{\ref {cquasi}}}+3\right\} \cap \left\{ d(x,\partial _{\ZZ ^{2}}M\dreg )>3\right\} \quad .\]
where $c_{\ref {cquasi}}$ is taken from lemma \ref{lem_no_quasi_loops}.
\label{random_hybrid} For every $0\leq k\leq \#Y$, define $H_{1,k}:=G(D_{k},N)$
to be a hybrid graph with $D_{k}$ a random subset of $Y$ of size
$k$. Let $a_{1,k}$ be the point of $H_{1,k}$ closest to $Ma$.
Let $S_{1,k}$ be a random walk on $H_{1,k}$ starting from $a_{1,k}$
and stopped on $B_{1,k}:=\partial _{H_{1,k}}M\dreg $. Let \[
p_{2,k}:=\PP (\LE (S_{1,k})\subset M\ereg )\]
(notice that $p_{1}=p_{2,0}$). \begin{sublem} \label{sublem_0_to_Y}
With the definitions above \begin{equation}
|p_{2,k}-p_{2,k+1}|\leq CM^{-c}\left(\frac{1}{k+1}+\frac{1}{\#Y-k}\right)+C\mu \quad .\label{step1}\end{equation}
\end{sublem} \begin{proof}[Subproof] We may couple $D_{k}$ and $D_{k+1}$
and assume that $D_{k}\subset D_{k+1}$ (one may either think about
$D_{k+1}$ as $D_{k}$ with a random point from $Y\setminus D_{k}$
added or about $D_{k}$ as $D_{k+1}$ with a random point removed).
Let $\{z\}:=D_{k+1}\setminus D_{k}$. Let $Z=z+[-1,2]^{2}$. We construct
$\LE (S_{1,k})$ as follows: 
\begin{itemize}
\item \label{page:gam123}let $b_{k}$ be a random point on $B_{1,k}$ chosen
with the hitting probabilities of $S_{1,k}$. 
\item Let $\check{S}_{k}$ be a random walk from $a_{1,k}$ to $B_{1,k}$
conditioned to hit $b_{k}$.
\item Let $\check{\gamma }_{k}$ be random simple path from $a_{1,k}$ to
$\left(\partial _{H_{1,k}}Z\right)\cup \{b_{k}\}$, which has the
same distribution as the segment of $\LE (\check{S}_{k})$ until $Z$,
or all of $\LE (\check{S}_{k})$ if $\LE (\check{S}_{k})\cap \partial Z=\emptyset $.
In particular, if $a_{1,k}\in Z$, then $\check{\gamma }_{k}=\{a_{1,k}\}$.
\item Let $c_{k}$ be the point where $\check{\gamma }_{k}$ hits $Z$ if
it does. If $a_{1,k}\in Z$ let $c_{k}=a_{1,k}$.
\item Let $\check{T}_{k}$ be a random walk on $H_{1,k}$ starting from
$b_{k}$ and conditioned to hit $B_{1,k}\cup \check{\gamma }_{k}$
in $c_{k}$, or $\emptyset $ if $\check{\gamma }_{k}$ never hits
$Z$.
\item Let $\gamma _{k}=\check{\gamma }_{k}\cup \LE (\check{T}_{k})$.
\end{itemize}
An easy application of lemma \ref{lem_condLE_sym} (symmetry of conditioned
loop-erased random walk) shows that $\gamma _{k}\sim \LE (S_{1,k})$.
Lemma \ref{lem_lerw_hit_local} shows that, if $N>CM^{C_{\ref {Cadmisspow}}}$,
\begin{eqnarray}
\sum _{b}|q_{k}^{1}-q_{k+1}^{1}| & \leq  & CN^{-1/3}\log N\log ^{2}M\leq C\mu \label{bk}\\
q_{k}^{1}(b) & := & \PP (b_{k}=b)\quad .\nonumber 
\end{eqnarray}
Next we use lemma \ref{lem_clerw_local} for the random walk on $H_{1,k}$
starting from $a_{1,k}$, stopped on $B_{1,k}$, and conditioned to
hit $b_{k}$ (notice that the condition $d(z,\linebreak [0]\partial _{\ZZ ^{2}}M\dreg )>3$
for all $z\in Y$ ensures the condition $B_{1,k}\cap (z+[-1,2]^{2})=\emptyset $
required by lemma \ref{lem_clerw_local}). This shows that \begin{eqnarray}
\sum _{\gamma }|q_{k}^{2}-q_{k+1}^{2}| & \leq  & C\mu \qquad \forall b\in B_{1,k}\label{gamma_k}\\
q_{k}^{2}(b,\gamma ) & := & \PP (\check{\gamma }_{k}=\gamma \, |\, b_{k}=b)\quad .\nonumber 
\end{eqnarray}
Thirdly, we use lemma \ref{lem_clerw_local}, this time for a random
walk starting from $b_{k}$, stopped on $B_{1,k}\cup \check{\gamma }_{k}$
and conditioned to hit $\check{c}_{k}$ to show that, when $\check{\gamma }_{k}'$
is the portion of $\LE (\check{T}_{k})$ up to $Z$,\begin{eqnarray}
|q_{k}^{3}-q_{k+1}^{3}| & \leq  & C\mu \qquad \forall b,\gamma \label{gamma_k_prime}\\
q_{k}^{3}(b,\gamma ) & := & \PP (\check{\gamma }_{k}'\subset M\ereg \, |\, b_{k}=b\; \wedge \; \gamma _{k}=\gamma )\quad .\nonumber 
\end{eqnarray}
Summing (\ref{bk}), (\ref{gamma_k}) and (\ref{gamma_k_prime}) gives\begin{eqnarray}
\lefteqn{|\PP ((\check{\gamma }_{k}\cup \check{\gamma }_{k}')\subset M\ereg )-\PP ((\check{\gamma }_{k+1}\cup \check{\gamma }_{k+1}')\subset M\ereg )|=} &  & \nonumber \\
 & \qquad  & =\left|\sum _{b,\gamma }q_{k}^{1}q_{k}^{2}q_{k}^{3}-q_{k+1}^{1}q_{k+1}^{2}q_{k+1}^{3}\right|\nonumber \\
 &  & \leq C\mu +\left|\sum _{b,\gamma }(q_{k}^{1}q_{k}^{2}-q_{k+1}^{1}q_{k+1}^{2})q_{k}^{3}\right|\nonumber \\
 &  & \leq C\mu \label{no_quasi_loop}
\end{eqnarray}
 In other words, we have proved that the probabilities (for $k$ and
$k+1$) that both segments of $\LE (S_{1,k})$, leading up to $Z$
and from $Z$ to $B_{1,k}$ to be in $\ereg $ are close. Thus the
only case we haven't covered is of $\LE (S_{1,k})$ doing a loop inside
$Z$. This would be a quasi-loop with $\epsilon \leq d(z,\partial _{H_{1,k}}Z)<3$.
Since $d(Z,M(\partial \ereg \setminus \partial \dreg ))>M^{1-c_{\ref {cquasi}}}$
we can denote\[
q_{k}^{4}(s):=\PP (\LE (S_{1,k})\in \QL (M^{1-c_{\ref {cquasi}}},3,s))\]
and then write (\ref{no_quasi_loop}) as\begin{equation}
|p_{2,k}-p_{2,k+1}|\leq C\mu +q_{k}^{4}(z)+q_{k+1}^{4}(z)\label{quasi_loop_left}\end{equation}

The estimate of (\ref{quasi_loop_left}) is where the random choice
of $z$ plays its part. Lemma \ref{lem_no_quasi_loops} gives us that
\[
\sum _{s\in Y\setminus D_{k}}q_{k}^{4}(s)\leq \sum _{s\in Y}q_{k}^{4}(s)\leq CM^{-c}\]
and since $z$ is chosen randomly from $Y\setminus D_{k}$ we have
\[
\EE q_{k}(z)\leq \frac{CM^{-c}}{\#(Y\setminus D_{k})}=\frac{CM^{-c}}{\#Y-k}\]
and therefore \[
\PP (\LE (S_{1,k})\in \QL (M^{1-c_{\ref {cquasi}}},3,z))=\EE q_{k}(z)\leq \frac{CM^{-c}}{\#Y-k}\quad .\]
For $q_{k+1}(z)$ we similarly have \[
\sum _{z\in D_{k+1}}q_{k+1}(z)\leq CM^{-c}\]
and since $z$ can also be thought of as being chosen randomly from
$D_{k+1}$ we get \[
\EE q_{k+1}(z)\leq \frac{CM^{-c}}{k+1}\]
and the sublemma is proved.\end{proof} 

\textbf{step 2}: Define \begin{eqnarray*}
p_{3} & := & \PP (\LE (S_{1,\#Y})\subset M\ereg _{2})\\
\ereg _{2} & := & \ereg +(2M^{-c_{\ref {cquasi}}}+6M^{-1})\DD \quad .
\end{eqnarray*}
Clearly this gives \begin{equation}
p_{3}>p_{2,\#Y}\label{step2}\end{equation}
and \[
d(M(\partial \ereg _{2}\setminus \partial \dreg ),Y')>M^{1-c_{\ref {cquasi}}}+3\]
where \[
Y':=X\setminus \left(Y\cup \{z\, :\, d(z,\partial _{\ZZ ^{2}}M\dreg )>3\right)\quad .\]

\textbf{step 3}: As in step 1, for every $0\leq k\leq \#Y'$, define
$H_{2,k}=G(D_{2,k},N)$ to be a hybrid graph with $D_{2,k}=Y\cup D_{k}'$
and $D_{k}'$ a random subset of $Y'$ of size $k$. Again, let $a_{2,k}$
be the point of $H_{2,k}$ closest to $Ma$, let $S_{2,k}$ be a random
walk on $H_{2,k}$ starting from $a_{2,k}$ and stopped on $\partial _{H_{2,k}}M\dreg $,
and let \[
p_{4,k}:=\PP (\LE (S_{2,k})\subset M\ereg _{2})\quad .\]
Again notice that $p_{4,0}=p_{3}$. \begin{sublem} With the definitions
above \begin{equation}
|p_{4,k}-p_{4,k+1}|\leq CM^{-c}\left(\frac{1}{k+1}+\frac{1}{\#Y'-k}\right)+C\mu \quad .\label{step3}\end{equation}
\end{sublem} The proof of this sublemma is identical to that of sublemma
\ref{sublem_0_to_Y} and we shall omit it.

\textbf{step 4:} \label{page:begstep4}Define $H_{3}=H_{2,\#Y'}$,
$a_{3}=a_{2,\#Y'}$ and $S_{3}$ a random walk on $H_{3}$ starting
from $a_{3}$ and stopped on $\partial _{H_{3}}M\dreg _{2}$ where\[
\dreg _{2}:=\CC \setminus \overline{(\CC \setminus \dreg )+5M^{-1}\DD }\quad .\]
If $a_{3}\not \in M\dreg _{2}$, let $S_{3}$ be the trivial path
$\{a_{3}\}$. Let\[
p_{5}:=\PP (\LE (S_{3})\subset M\ereg _{3})\]
where\[
\ereg _{3}:=\ereg _{2}+M^{-c_{\ref {cquasi}}}\DD \]
\begin{sublem} \label{sublem:step4}For $M$ sufficiently large\begin{equation}
p_{5}>p_{4,\#Y'}-CM^{-c}\quad .\label{eq:step4}\end{equation}
\end{sublem} \begin{proof}[Subproof] $\CC \setminus \dreg $ has
a finite number of connected components, $\{T_{i}\}$. The quantity
that interests us is \[
\tau :=\min _{T_{i}}\left\{ \diam T_{i}\right\} \quad .\]
Now the walks $S_{2,\#Y'}$ and $S_{3}$ are walks on the same graph
stopped at $\partial M\dreg $ and $\partial M\dreg _{2}$ respectively.
Therefore if we define $t_{1}$ and $t_{2}$ to be the stopping times
of $S_{3}$ on $\partial M\dreg $ and $\partial M\dreg _{2}$ (define
$t_{2}=0$ if $a\not \in M\dreg _{2}$) then the question reduces
to an estimate of \begin{equation}
\PP (\LE (S_{3}([0,t_{2}]))\not \subset M\ereg _{3}\; \wedge \; \LE (S_{3}([0,t_{1}]))\subset M\ereg _{2})\quad .\label{eq:St2_large_St1_small}\end{equation}
Let $T$ be the graph-connected-component of $H_{3}\setminus (M\dreg )^{\circ }$
closest to $S_{3}(t_{2})$ --- the definition of $\dreg _{2}$ gives
that $d(T,S_{3}(t_{2}))\leq 6$. It's easy to see that $\diam T\geq \tau M-3$,
and then get from lemma \ref{lem_hybrid_Kesten} that \begin{equation}
\PP (S_{3}([t_{2},t_{1}])\text {\sh exits\sh }S_{3}(t_{2})+\epsilon \DD )\leq C\left(\min (\tau M-3,\epsilon )\right)^{-c}\quad \forall \epsilon \geq 1.\label{eq:r_large}\end{equation}
On the other hand, if $S_{3}([t_{2},t_{1}])\subset S_{3}(t_{2})+\epsilon \DD \subset S_{3}(t_{1})+2\epsilon \DD $
and in addition the event of (\ref{eq:St2_large_St1_small}) hold
then we can conclude that $\LE (S_{3}[0,t_{2}])\in \QL (M^{1-c_{\ref {cquasi}}},2\epsilon ,S_{3}(t_{1}))$,
and lemma \ref{lem_no_quasi_loops} gives the bound \begin{equation}
\PP \left(\LE (S_{3}[0,t_{2}])\in \QL \left(M^{1-c_{\ref {cquasi}}},2\epsilon ,S_{3}(t_{1})\right)\right)\leq CM^{-c}\epsilon ^{C}\quad .\label{eq:r_small}\end{equation}
We choose $\epsilon =M^{c}$ with $c$ sufficiently small and combine
(\ref{eq:r_large}) and (\ref{eq:r_small}) to get the required estimate
of (\ref{eq:St2_large_St1_small}) which holds whenever $\epsilon \leq \tau M-3$
or equivalently \begin{equation}
M>\tau ^{-c}+C\quad .\qedhere \label{eq:M_Dgeo}\end{equation}
\end{proof} 

\textbf{step 5}: Define $H_{4}=G([-M,M]^{2}\cap \ZZ ^{2},N)$, $a_{4}$
the point of $H_{4}$ closest to $a$ and $S_{4}$ a random walk on
$H_{4}$ starting from $a_{4}$ and stopped on $\partial _{H_{4}}M\dreg $.
Let\[
p_{6}:=\PP (\LE (S_{4})\subset M\ereg _{4})\]
where\[
\ereg _{4}:=\ereg _{3}+M^{-1/2}\DD \]
\begin{sublem} With the definitions above\begin{equation}
p_{6}>p_{5}-CM^{-c}\quad .\label{eq:step5}\end{equation}
\end{sublem} 
\begin{proof}
[Subproof]As in sublemma \ref{sublem:step4}, we need to show that\[
\PP (\LE (S_{4}([0,t_{1}])\not \subset M\ereg _{4}\; \wedge \; \LE (S_{4}([0,t_{2}])\subset M\ereg _{3})\leq CM^{-c}\]
with the same $t_{1}$ and $t_{2}$. Unlike in sublemma \ref{sublem:step4},
this requires no recourse to lemma \ref{lem_no_quasi_loops} but rather
follows directly from lemma \ref{lem_hybrid_Kesten} since this event
implies that $S_{4}[t_{2},t_{1}]\not \subset S_{4}(t_{2})+M^{1/2}\DD $
whose probability can be bounded by \[
C\left(\min (\tau M-3,M^{1/2})\right)^{-c}\]
and if (\ref{eq:M_Dgeo}) is fulfilled then this is $\leq CM^{-c}$.\label{page:endstep5}
\end{proof}
\textbf{final step}: At this point our environment is no longer hybrid%
\footnote{Actually, it was already true in step 5.%
} --- in effect $H_{4}\cap \dreg \equiv (\frac{1}{2}M\delta )\ZZ ^{2}\cap \dreg $.
Thus we can return to the notations of $G_{i}$, $R_{i}$ etc. and
get \[
p_{6}=\PP (\LE (R_{2})\subset \ereg _{4})\]
Summing up (\ref{step1}), (\ref{step2}), (\ref{step3}), (\ref{eq:step4})
and (\ref{eq:step5}) we get \newc{c:Mtodel} \begin{equation}
p_{6}>p_{1}-CM^{-c}\log M+CM^{2}\mu \quad .\label{just_choose_M}\end{equation}
The only thing left now is to choose $M$. The following conditions
must be met:
\begin{enumerate}
\item $CM^{2}\mu \leq CM^{-c}$;
\item $N>CM^{C_{\ref {Cadmisspow}}}$ --- this will also give that $H_{1,k}$
and $H_{2,k}$ are admissible;
\item $N\in \NN $;
\item $M>\tau ^{-1/2}+C$ (that's (\ref{eq:M_Dgeo}) \vpageref[above]{eq:M_Dgeo}).
\end{enumerate}
For some $c_{\ref {c:Mtodel}}$ sufficiently small, if we choose $M\approx \delta ^{-c_{\ref {c:Mtodel}}}$
then we will have $N\approx M^{\frac{1-c_{\ref {c:Mtodel}}}{c_{\ref {c:Mtodel}}}}$
and therefore (say take $c_{\ref {c:Mtodel}}<\frac{1}{7}$) that $M^{2}\mu <CM^{-c}$.
Requirement 2 will also follow if $c_{\ref {c:Mtodel}}$ is sufficiently
small --- this depends on the constant $C_{\ref {Cadmisspow}}$ that
appears in lemma \ref{lemma_a_wired}. Since $C_{\ref {Cadmisspow}}$
can be chosen to be any value $>2$ then the restriction on our $c_{\ref {c:Mtodel}}$
is in effect only the weaker $c_{\ref {c:Mtodel}}<\frac{1}{3}$. To
fulfill condition 4, we need some assumption on $\delta $: $\delta <\delta _{0}(\dreg )=c\tau ^{C}$
will be enough. Clearly condition 3 is no obstacle. Plugging this
into (\ref{just_choose_M}) will give \[
p_{6}>p_{1}-CM^{-c}>p_{1}-C\delta ^{c}\quad .\]
On the other hand, for an appropriate $C_{\ref {Cmain}}$ and $c_{\ref {cmain}}$,
$\ereg _{4}\subset \ereg +C_{\ref {Cmain}}\delta ^{c_{\ref {cmain}}}\DD $.
This finishes (\ref{R2_geq_R1}).

The proof of (\ref{R1_geq_R2}) is identical, with $D$ and $([-M,M]\cap \ZZ ^{2})\setminus D$
replaced everywhere. Thus ends the main lemma.
\end{proof}
\renewcommand{\thesublem}{\arabic{lem}.\arabic{sublem}} \begin{rmks}

\begin{enumerate}
\item The requirement from $\dreg $ to be a polygon was rather excessive.
In effect we used it only in steps 4 and 5 to show $\tau >0$. Therefore
the main lemma holds, for example, for any bounded domain $\dreg $
with no punctures (here we mean punctures in the sense of connected
components of $\CC \setminus \dreg $ with only one point, but not
necessarily isolated). It is not difficult to see that punctures in
$\dreg $ would require to reformulate the main lemma so as to take
into consideration the distance between $a$ and the nearest puncture.
See also the example on page \pageref{page:exam_punct} for the problems
punctures could bring about.
\item The division into $z$'s close to $\partial M\dreg $ and far from
it is not really necessary --- it is possible to extend lemmas \ref{lem_lerw_hit_local}
and \ref{lem_clerw_local} to work when $B\cap Z\neq \emptyset $
and thus save steps 4 and 5 in the main lemma. However, with this
extension the formulation of lemmas \ref{lem_lerw_hit_local} and
\ref{lem_clerw_local} is very awkward. We would need two $B_{i}$
which are {}``almost similar'', two $b_{i}$'s, and make provisions
for the cases when $\kappa _{b_{1}}\neq \kappa _{b_{2}}$ since the
probability to hit $b_{i}$ depends on $\kappa _{b_{i}}$ (see (\ref{def_kappab})).
The proofs (especially that of lemma \ref{lem_lerw_hit_local}) would
also suffer from a canworm of geometric issues.
\item An alternative to the use of random hybrid graphs, is to randomize
the starting point $a$. This would give similar results (especially
with results of the next section).
\item I am happy to promise to my readers that this is the last time the
term {}``hybrid graph'' is mentioned in this paper. Or, to be more
precise, we will still refer to some lemmas formulated using hybrid
graphs --- particularly to the ubiquitous lemma \ref{lem_hybrid_Kesten}
--- but only for the non-hybrid case i.e.~$D=\emptyset $.
\end{enumerate}
\end{rmks}

\subsection{\label{sec:cont}Continuity }

In this section we prove some simple estimates that show that the
probability of a loop-erased random walk to be in a set is continuous
in the point of departure, the set and the environment.

\begin{lem}
\label{lem_cont_a} Let $\ereg $ and $\dreg $ be open sets. Let
$v,w\in \ereg \cap \dreg \cap \ZZ ^{2}$. Let $R^{x}$ be a random
walk on $\ZZ ^{2}$ started from $x$ and stopped on $\partial _{\ZZ ^{2}}\dreg $.
Then \[
|\PP (\LE (R^{v})\subset \ereg )-\PP (\LE (R^{w})\subset \ereg )|\leq C\left(\frac{d(v,\partial \ereg \cup \partial \dreg )}{d(v,w)}\right)^{-c}\]

\end{lem}
\begin{proof}
Denote $\mu =d(v,\partial \ereg \cup \partial \dreg )/d(v,w)$. Let
$S^{w}$ be a random walk started from $w$ and stopped on $\LE (R^{v})\cup \partial G$.
Lemma \ref{lem_hybrid_Kesten} says that the probability of $S^{w}$
to hit $R^{v}$ before exiting $\ereg \cap \dreg $ is $\leq C\mu ^{-c}$.
But since Wilson's algorithm says that $\LE (R^{w})$ has the same
distribution as $\LE (S^{w})$ unioned with the segment of $\LE (R^{v})$
from $S^{w}\cap \LE (R^{v})$ to $\partial G$, the lemma is finished. 
\end{proof}
\begin{defn*}
For $\ereg ,\dreg $ open and for $r>0$ we define $X_{1}(r;\dreg )\subset \ZZ ^{2}\setminus \dreg ^{\circ }$
to be the union of all graph-connected-components of $\ZZ ^{2}\setminus \dreg ^{\circ }$
of diameter $<r$. Next we define \begin{equation}
X_{2}(r;\ereg ,\dreg ):=\{x\in \partial _{\ZZ ^{2}}\dreg \, :\, d(x,\partial \ereg )<r\}\label{eq:defX2_forrho2}\end{equation}
and thirdly $X_{3}:=X_{1}\cup X_{2}$. Next, for $a\in \dreg \cap \ZZ ^{2}$
and for $i=1,2,3$, define \[
\rho _{i}(r;a,\ereg ,\dreg ):=q(a,X_{i},\partial \dreg ,\ZZ ^{2})\quad .\]
and\[
\rho _{i}(r,\delta ;a,\ereg ,\dreg ):=\rho _{i}(\delta ^{-1}r;\left\lfloor \delta ^{-1}a\right\rfloor ,\delta ^{-1}\ereg ,\delta ^{-1}\dreg )\quad .\]

\end{defn*}
The {}``good'' sets (or rather triplets $a$, $\ereg $, $\dreg $)
are the ones satisfying \[
\lim _{(r,\delta )\rightarrow (0,0)}\rho _{3}(r,\delta )=0\quad .\]
There are counter example, though. In the example on page \pageref{page:exam_punct},
with $\ereg =\left]-2,2\right[^{2}$, we have $\limsup \rho _{1}>0$.
It is also possible to construct non-trivial examples where the culprit
is $\rho _{2}$.

\begin{lem}
\label{lem_cond_D}Let $0<s<r$. Let $\dreg _{1}$, $\dreg _{2}$
and $\ereg $ be open sets and assume that $\dreg _{1}$ and $\dreg _{2}$
are similar in the following sense:\begin{equation}
d(v,\partial \dreg _{1}\setminus X_{3}(r;\ereg ,\dreg _{1}))\leq s\quad \forall v\in \partial \dreg _{2}\setminus X_{3}(r+s;\ereg ,\dreg _{2})\label{eq:X3simcond}\end{equation}
and similarly with $\dreg _{1}$ and $\dreg _{2}$ replaced. Let $a\in \dreg _{1}\cap \dreg _{2}\cap \ereg \cap \ZZ ^{2}$.
Let $R_{i}$ be random walks on $\ZZ ^{2}$ started from $a$ and
stopped on $\partial _{\ZZ ^{2}}\dreg _{i}$. Then\begin{eqnarray*}
|p_{1}-p_{2}| & \leq  & \rho _{3}(r+s;a,\ereg ,\dreg _{1})+\rho _{3}(r+s;a,\ereg ,\dreg _{2})+C\left(\frac{r}{s}\right)^{-c}\\
p_{i} & := & \PP (\LE (R_{i})\subset \ereg )\quad .
\end{eqnarray*}
If $\delta \neq 1$ this is true if (\ref{eq:X3simcond}) holds with
a $\delta $ version of $X_{3}$: $X_{3}(r,\delta ;\ereg ,\dreg ):=\delta X_{3}(\delta ^{-1}r;\delta ^{-1}\ereg ,\delta ^{-1}\dreg )$.
\end{lem}
\begin{proof}
This follows from lemmas \ref{lem_hybrid_Kesten} and \ref{lem_no_quasi_loops}
--- see steps 4 and 5 of the main lemma \vpagerefrange{page:begstep4}{page:endstep5}
for a more detailed version of this argument.
\end{proof}
\begin{lem}
\label{lem_cont_E_move}Let $\ereg ,\dreg $ be open sets. Let $v\in \dreg \cap \ZZ ^{2}$.
Let $R_{\delta }$ be a random walk on $\delta \ZZ ^{2}$ started
from $v$ and stopped on $\partial \dreg $. Then \[
\limsup _{\substack{ (\tau ,\delta )\rightarrow (0,0)\brk \tau \in \delta \ZZ ^{2}}
}|\PP (\LE (R_{\delta })\subset \ereg )-\PP (\LE (R_{\delta })\subset \ereg +\tau )|\leq 3\limsup _{(r,\delta )\rightarrow (0,0)}\rho _{3}(r,\delta ;a,\ereg ,\dreg )\]

\end{lem}
\begin{proof}
Let $R_{\delta }'$ be a random walk stopped on $\partial \dreg -\tau $.
We use lemma \ref{lem_cond_D} with $s=|\tau |$ and some $r>|\tau |$
to get\begin{eqnarray}
\lefteqn{|\PP (\LE (R_{\delta })\subset \ereg )-\PP (\LE (R_{\delta }')\subset \ereg )|\leq } &  & \nonumber \\
 & \qquad  & \leq \rho _{3}(r+|\tau |,\delta ;\dreg )+\rho _{3}(r+|\tau |,\delta ;\dreg -\tau )+C\left(\frac{r}{|\tau |}\right)^{-c}\quad .\label{eq:Rdel_Rdelp}
\end{eqnarray}
Now to estimate $\rho _{3}(\dreg -\tau )$, we write $\rho _{3}\leq \rho _{1}+\rho _{2}$
and an argument like lemma \ref{lem_cont_a} shows\begin{equation}
\rho _{1}(r,\delta ;a,\dreg -\tau )=\rho _{1}(r,\delta ;a+\tau ,\dreg )\leq \rho _{1}(r,\delta ;a,\dreg )+C\left(\frac{d(a,\partial \dreg )}{\tau }\right)^{-c}\label{eq:rho1_mov}\end{equation}
and similarly\begin{eqnarray}
\rho _{2}(r,\delta ;a,\dreg -\tau ,\ereg ) & \leq  & \rho _{2}(r+|\tau |,\delta ;a+\tau ,\dreg ,\ereg )\leq \nonumber \\
 & \leq  & \rho _{2}(r+|\tau |,\delta ;a,\dreg ,\ereg )+C\left(\frac{d(a,\partial \dreg \cup \partial \ereg )}{\tau }\right)^{-c}\quad .\label{eq:rho2_mov}
\end{eqnarray}
 Next we define $R_{\delta }''$ to be a random walk starting from
$a-\tau $ stopped on $\dreg -\tau $ and again use lemma \ref{lem_cont_a}
to get \begin{equation}
|\PP (\LE (R_{\delta }')\subset \ereg )-\PP (\LE (R_{\delta }'')\subset \ereg )|\leq C\left(\frac{d(a,\partial \dreg \cup \partial \ereg )}{|\tau |}\right)^{-c}\quad .\label{eq:Rdelp_Rdelpp}\end{equation}
Summing (\ref{eq:Rdel_Rdelp}), (\ref{eq:rho1_mov}), (\ref{eq:rho2_mov})
and (\ref{eq:Rdelp_Rdelpp}) and estimating $\rho _{1}$,$\rho _{2}\leq \rho _{3}$
we get \begin{eqnarray}
\lefteqn{|\PP (\LE (R_{\delta })\subset \ereg )-\PP (\LE (R_{\delta })\subset \ereg +\tau )|=} &  & \label{eq:Rdel_Rdelpp}\\
 & \qquad \qquad  & =|\PP (\LE (R_{\delta })\subset \ereg )-\PP (\LE (R_{\delta }'')\subset \ereg )|\leq \nonumber \\
 &  & \leq 3\rho _{3}(r+2|\tau |,\delta ;a,\ereg ,\dreg )+C\left(\frac{d(a,\partial \dreg \cup \partial \ereg )}{|\tau |}\right)^{-c}+C\left(\frac{r}{|\tau |}\right)^{-c}\nonumber 
\end{eqnarray}
and choosing $r=\sqrt{|\tau |}$ will make the two summands on the
right of (\ref{eq:Rdel_Rdelpp}) converge to $0$ when $(\tau ,\delta )\rightarrow (0,0)$.
\end{proof}
\begin{lem}
\label{lem_cont_E}Let $a\in \ereg \cap \dreg $ where $\ereg $ is
a polygon, and assume \begin{equation}
\lim _{(r,\delta )\rightarrow (0,0)}\rho _{3}(r,\delta ;a,\ereg ,\dreg )=0\quad .\label{always}\end{equation}
Let $R_{\delta }$ be a random walk on $G:=\delta \ZZ ^{2}$ starting
from $a$ and stopped on $\partial _{G}\dreg $. Then

\begin{equation}
\lim _{(\epsilon ,\delta )\rightarrow (0,0)}|\PP (\LE (R_{\delta })\subset \ereg )-\PP (\LE (R_{\delta })\subset \ereg +\epsilon \DD )|=0\quad .\label{eq:lem_cont_E}\end{equation}

\end{lem}
\begin{proof}
For $A,B$ satisfying $A\cap \dreg \subset B\cap \dreg $ we denote
\[
F(A,B):=\PP (\emptyset \neq (\LE (R_{\delta })\cap B)\subset A)\]
so (\ref{eq:lem_cont_E}) is equivalent to \[
\lim _{(\epsilon ,\delta )\rightarrow (0,0)}F\left((\ereg +\epsilon \DD )\setminus \ereg ,\dreg \setminus \ereg \right)=0\quad .\]
Let $\left\{ \xi _{i}\right\} _{i=1}^{n}$ be the vertices of $\ereg $
inside $\dreg $, let $\left\{ I_{i}\right\} _{i=1}^{n+1}$ be the
segments and let $I_{i}^{\bot }$ be orthogonal segments oriented
to the exterior of $\ereg $. Denote $I_{i}^{\bot }=[0,\nu _{i}[$,
$|\nu _{i}|=1$. For each $\xi _{i}$ we use a simple hitting probability
estimate to get \begin{eqnarray*}
\PP (\LE (R_{\delta })\cap (\xi _{i}+\epsilon \DD )\neq \emptyset ) & \leq  & \PP (R_{\delta }\cap (\xi _{i}+\epsilon \DD )\neq \emptyset )\\
 & \leq  & C\frac{\log (d(a,\xi _{i})/\diam \dreg )}{\log (\epsilon /\diam \dreg )}
\end{eqnarray*}
which converges to zero (in effect, much more precise estimates are
known --- see e.g. \cite{K00}). Let therefore $\epsilon _{1}(\mu )$
and $\delta _{1}(\mu )$ satisfy that\[
\sum _{i}\PP (\LE (R_{\delta })\cap (\xi _{i}+\epsilon \DD )\neq \emptyset )<\mu \quad \forall \delta <\delta _{1}(\mu ),\: \epsilon <\epsilon _{1}(\mu )\]
which gives us in $F$ notation\[
F((\ereg +\epsilon \DD )\setminus \ereg ,\dreg \setminus \ereg )\leq \mu +F\left(\bigcup P_{i}(\epsilon ),\dreg \setminus \ereg \right)\]
where $P_{i}(\epsilon ):=\left(I_{i}+\epsilon I_{i}^{\bot }\right)\setminus \bigcup (\xi _{i}+\epsilon _{1}(\mu )\DD )$.
Clearly, the event $\emptyset \neq \LE (R_{\delta })\cap (\dreg \setminus \ereg )\subset \bigcup P_{i}(\epsilon )$
implies that for at least one $i$, $\emptyset \neq \LE (R_{\delta })\cap P_{i}(\epsilon )$.
If in addition the different $P_{i}(\epsilon )$ are disjoint, which
happens for $\epsilon <\epsilon _{2}(\mu )$, then we have \[
\leq \mu +\sum _{i}F\bigg (P_{i}(\epsilon ),\dreg \setminus \Big (\ereg \cup \bigcup _{j\neq i}P_{j}(\epsilon )\Big )\bigg )\quad .\]

The next step are some $F$ calculations, based on the following easy
inequalities which we call the monotonicity of $F$:\begin{equation}
\begin{array}{rcl}
 F(A,B) & \leq  & F(A\cup C,B\cup C)\\
 F(A,B) & \geq  & F(A,B\cup C)\quad .\end{array}
\label{F_monotone}\end{equation}
For $r=1,2,3$ we define\[
J_{i}^{r}:=I_{i}\setminus \bigcup _{j}\left(\xi _{j}+\quarter r\epsilon _{1}(\mu )\DD \right)\quad .\]
 For a suitable $\epsilon _{3}(\mu )$ we shall get $P_{i}(\epsilon )\subset J_{i}^{3}+\epsilon I_{i}^{\bot }$
whenever $\epsilon <\epsilon _{3}(\mu )$, and $(P_{i}(\epsilon )\cup \ereg )\cap (J_{j}^{1}+dI_{j}^{\bot })=\emptyset $
whenever $i\neq j$ and $\epsilon <d<\epsilon _{3}(\mu )$. This allows
to write \[
F\bigg (P_{i}(\epsilon ),\dreg \setminus \Big (\ereg \cup \bigcup _{j\neq i}P_{j}(\epsilon )\Big )\bigg )\leq F(P_{i}(\epsilon ),J_{i}^{1}+dI_{i}^{\bot })\leq F(J_{i}^{3}+\epsilon I_{i}^{\bot },J_{i}^{1}+dI_{i}^{\bot })\]
From this point on we shall drop the notation $i$ from $I_{i}^{\bot }$,
$J_{i}^{r}$, and $\nu _{i}$. Now, lemma \ref{lem_cont_E_move} with
(\ref{always}) and the obvious \[
\rho _{3}(\epsilon ;J^{1}+dI^{\bot })\leq \rho _{3}(\epsilon +d;\ereg )\]
give that for any $\epsilon <d<\epsilon _{4}(\mu )$, $\delta <\delta _{2}(\mu )$
and for any $\lambda \in \delta \ZZ ^{2}$, $|\lambda |<\lambda _{1}(\mu )$,
\begin{equation}
|F(J^{3}+\epsilon I^{\bot },J^{1}+dI^{\bot })-F(J^{3}+\lambda +\epsilon I^{\bot },J^{1}+\lambda +dI^{\bot })|\leq \mu \quad .\label{J_move}\end{equation}
Next, if in addition $\epsilon <\frac{d}{2m},\frac{\lambda _{1}(\mu )}{2m}$
and $\delta <\delta _{3}(m,\epsilon )$ we can pick $\lambda _{1},...,\lambda _{m-1}\in \delta \ZZ ^{2}$,
$|\lambda _{i}|\leq \lambda _{1}(\mu )$ that will satisfy%
\begin{figure}
[b]\input{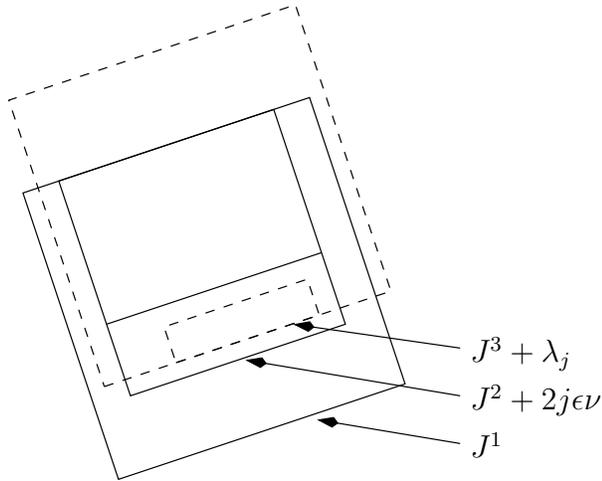}

\caption{\label{fig_diskette}Segments shifted by multiples of $\epsilon \nu $
are in solid lines; segments shifted by $\lambda _{j}$ are in dashed
lines.}
\end{figure}
\begin{eqnarray*}
J^{3}+\lambda _{j}+\epsilon I^{\bot } & \subset  & J^{2}+2j\epsilon \nu +2\epsilon I^{\bot }\\
J^{1}+\lambda _{j}+dI^{\bot } & \supset  & J^{2}+2(j+1)\epsilon \nu +2(m-j-1)\epsilon I^{\bot }
\end{eqnarray*}
 (just take $\lambda _{j}$ similar to $(2j+\frac{1}{2})\epsilon \nu $).
See figure \ref{fig_diskette} \vpageref[below]{fig_diskette}. From
these and (\ref{F_monotone}) we get\begin{eqnarray}
\lefteqn{F(J^{3}+\lambda _{j}+\epsilon I^{\bot },J^{1}+\lambda _{j}+dI^{\bot })\leq } &  & \nonumber \\
 & \qquad  & \leq F(J^{3}+\lambda _{j}+\epsilon I^{\bot },(J^{1}+\lambda _{j}+dI^{\bot })\cap (J^{2}+2j\epsilon \nu +2(m-j)\epsilon I^{\bot })\nonumber \\
 &  & \leq F(J^{2}+2j\epsilon \nu +2\epsilon I^{\bot },J^{2}+2j\epsilon \nu +2(m-j)\epsilon I^{\bot })\quad .\label{J_resize}
\end{eqnarray}
The same holds for $j=0$ with $\lambda _{0}=0$. But clearly \begin{eqnarray}
\lefteqn{F(J^{2}+2m\epsilon I^{\bot },J^{2}+dI^{\bot })=} &  & \nonumber \\
 & \qquad  & \sum _{j=0}^{m}F(J^{2}+2j\epsilon \nu +2\epsilon I^{\bot },J^{2}+2j\epsilon \nu +(d-2j\epsilon )I^{\bot })\label{F_is_sum_F}
\end{eqnarray}
since the event $\emptyset \neq \LE (R_{\delta })\cap J^{2}+dI^{\bot }\subset J^{2}+2m\epsilon I^{\bot }$
whose probability is measured on the left can be divided according
to $\max \{j:R_{\delta }\subset J^{2}+2j\epsilon I^{\bot }\}$. Summing
(\ref{J_move}), (\ref{J_resize}) and (\ref{F_is_sum_F}) we get
\[
1\geq F(J^{2}+2m\epsilon I^{\bot },J^{2}+dI^{\bot })\geq mF(J^{3}+\epsilon I^{\bot },J^{1}+dI^{\bot })-m\mu \quad .\]
This finishes the lemma --- we pick $m$ large, then pick $\mu <\frac{1}{m}$,
then $d<\epsilon _{4}(\mu ),\lambda _{1}(\mu )$, then $\epsilon _{\min }:=\min \epsilon _{1}(\mu ),\epsilon _{2}(\mu ),\epsilon _{3}(\mu ),\frac{d}{2m}$
and for every $\delta <\delta _{1}(\mu ),\linebreak [0]\delta _{2}(\mu ),\linebreak [0]\delta _{3}(m,\epsilon _{\min })$
we get \[
F(J^{1}+\epsilon _{\min }I^{\bot },J^{3}+dI^{\bot })\leq \frac{C}{m}\]
and since $F(J^{1}+\epsilon I^{\bot },J^{3}+dI^{\bot })$ is decreasing
in $\epsilon $ and since $F((\ereg +\epsilon \DD )\setminus \ereg ,\dreg \setminus \ereg )\leq $
a finite sum of those, we are done. 
\end{proof}
\begin{lem}
\label{lem:limit}Let $a\in \ereg \cap \dreg $ where $\ereg $ and
$\dreg $ are polygons. Let $R_{n}$ be a random walk on $G:=2^{-n}\ZZ $
starting from $a$ and stopped on $\partial _{G}\dreg $. Then $\PP (\LE (R_{n})\subset \ereg )$
converges to a limit as $n\rightarrow \infty $.
\end{lem}
\begin{proof}
Use the main lemma (shrink $\dreg $ to fit $[-1,1]^{2}$ if necessary)
repeatedly to get $\forall n,m>N_{0}$ \[
\PP (\LE (R_{n})\subset \ereg )\leq \PP (\LE (R_{m})\subset \ereg +C2^{-cN_{0}}\DD )+C2^{-cN_{0}}\quad .\]
 Since it is obvious that for $\dreg $ and $\ereg $ polygons $\rho (\epsilon )\rightarrow 0$
(if $\partial \dreg \cap \partial \ereg $ contains a segment, just
extend $\ereg $ a bit so as to make $\partial \dreg \cap \partial \ereg $
finite), lemma \ref{lem_cont_E} gives \[
\PP (\LE (R_{n})\subset \ereg +C2^{-cN_{0}}\DD )\leq \PP (\LE (R_{n})\subset \ereg )+\mu (N_{0})\]
 where $\mu (N_{0})\rightarrow 0$, which gives\[
|\PP (\LE (R_{n})\subset \ereg )-\PP (\LE (R_{m})\subset \ereg )|\leq C2^{-cN_{0}}+\mu (N_{0})\quad .\qedhere \]

\end{proof}
\begin{rem*}
My intuition would have been that lemma \ref{lem:limit} --- while
probably being worthy of proof in its own right --- wouldn't be necessary
for the proof of the theorem since the weak limit is insensitive to
small inflations. However, I was not able to surmount certain technical
difficulties in translating this intuition into a proper proof.
\end{rem*}

\section{\label{chap_limit}The limit process}

In this section we shall conclude the theorem from lemma \ref{lem:limit}.
This is a standard limit process, so we shall explain it only briefly.
We start with an example that will justify our choice of acceptable
$\dreg $'s. 

\begin{example*}
\label{page:exam_punct}An open set $\dreg \subset \CC $ such that
the loop-erased random walk (and even the regular random walk) from
a point $a$ to $\partial _{2^{-n}\ZZ ^{2}}\dreg $ does not converge
as $n\rightarrow \infty $ in $\Meas (\Haus (\dreg ))$.
\begin{itemize}
\item Let $E'\subset \left]-\frac{1}{3},\frac{1}{3}\right[$ be a closed
set with $mE'\geq \frac{1}{2}$ and $E'\cap \QQ =\emptyset $, $\QQ $
being the rationales so $E'$ is a Cantor-like set. Let \[
E=\left(E'\pm i\tsh {\frac{1}{3}}\right)\cup \left(iE'\pm \tsh {\frac{1}{3}}\right)\quad .\]

\item Let $n_{k}$ be a sequence converging to $\infty $ sufficiently fast
(we shall specify them later).
\item Let $r_{k}\in 2^{-n_{k}}\ZZ \setminus 2^{-n_{k}+1}\ZZ $ satisfy $|r_{k}-\frac{1}{3}|\leq 2^{-n_{k}+1}$.
\item Let\[
P_{k}':=\left(2^{-n_{k}}\ZZ \setminus 2^{-n_{k}+1}\ZZ \right)\cap \left\{ x\, :\, d(x,E')<\frac{1}{k}\right\} \]
\[
P_{k}:=(P_{k}'\pm ir_{k})\cup (iP_{k}'\pm r_{k})\]

\item Let \[
\dreg =[-1,1]^{2}\setminus \left(E\cup \bigcup _{k=1}^{\infty }P_{k}\right)\]

\end{itemize}
It is easy to see that $\dreg $ is an open set.
\end{example*}
\begin{lem}
With an appropriate choice of $n_{k}$, the probabilities \[
q(0,\partial [-1,1]^{2},\partial \dreg ,2^{-n}\ZZ ^{2})\]
do not converge.
\end{lem}
This of course implies that the distributions of the regular random
walk do not converge in $\Meas (C([0,\infty [\rightarrow \dreg ))$,
the distributions of the loop-erased random walk do not converge in
$\Meas (\Haus (\dreg ))$ and would probably exclude convergence in
any reasonable topology.

\begin{proof}
Denote $\dreg _{k}=[0,1]^{2}\setminus \bigcup _{j=1}^{k}P_{j}$. Because
$P_{k}'$, $r_{k}$ and $E'$ all avoid $2^{-n}\ZZ $ for all $n<n_{1}$,
and moreover, avoid all edges of this graph (when viewed as line segments
in $\CC $), we get that $R_{n}$ is identical to a walk on $2^{-n}\ZZ ^{2}\cap \dreg _{0}\equiv 2^{-n}\ZZ ^{2}\cap [-1,1]^{2}$.
At $n_{1}$ we have that $\partial \dreg \subset \partial [-r_{k},r_{k}]^{2}$
satisfies that $\#\partial \dreg \geq c\#\partial [-r_{k},r_{k}]^{2}$
and in particular has a hitting probability $\geq c$. For $n_{1}\leq n<n_{2}$,
however, $\partial \dreg =\partial \dreg _{1}$ and since it is just
a finite set of points, for $n_{2}$ sufficiently large the hitting
probability of $\partial \dreg _{1}$ can be made as small as desired.
For $n=n_{2}$ we have $\partial \dreg =\partial \dreg _{2}$ and
again has hitting probability $\geq c$, etc. To sum it all up, if
$n_{k}$ increases sufficiently fast, we have that the probability
of $R_{n}$ to hit $\partial [-1,1]^{2}$ fluctuates between $1$
and $c$. This $c$ can be chosen arbitrarily close to $0$ merely
by changing $E'$.
\end{proof}
In view of this example, we must somehow restrict the $\dreg $'s
we talk about. We shall examine the class of bounded finitely-connected
open sets.

\begin{lem}
\label{lem_limit}The conclusion of lemma \ref{lem:limit} holds for
any open, bounded and finitely connected $\dreg $ if \[
\lim _{(\epsilon ,n)\rightarrow (0,\infty )}\rho _{3}(\epsilon ,2^{-n};a,\dreg ,\ereg )\rightarrow 0\quad .\]

\end{lem}
\begin{proof}
Let $\left\{ H_{i}\right\} _{i=0}^{m}$ be the connected components
of $\CC \setminus \dreg $, $H_{0}$ being the unbounded one. Let
$\mu >0$ be some parameter. Let $P_{i}$ be simply connected polygons
with $d_{\Haus }(\partial H_{i},\partial P_{i})\leq \mu $. Assume
$\mu $ is sufficiently small as to make the $\partial P_{i}$ pairwise
disjoint, and define $P=P_{0}\setminus \cup _{i>0}P_{i}$. Our goal
is to use lemma \ref{lem_cond_D} for $\dreg $ and $P$, and we need
to estimate the effects of discretization. This is easy to do, since
for any $X$ simply connected, \[
\partial _{\delta \ZZ ^{2}}X\neq \emptyset \Rightarrow d_{\Haus }(\partial X,\partial _{\delta \ZZ ^{2}}X)\leq \tsh {\frac{1}{\sqrt{2}}}\delta \]
and if $Y=\left(X\cap \delta \ZZ ^{2}\right)\cup \partial _{\delta \ZZ ^{2}}X$
then \[
\left|\diam X-\diam Y\right|\leq \sqrt{2}\delta \quad .\]
These two imply that, when $2^{-n}<\mu $, the requirements of lemma
\ref{lem_cond_D} will be satisfied when $s>C\mu $ and as a result
we will get, \begin{eqnarray*}
\lefteqn{|\PP (\LE (R_{n})\subset \ereg )-\PP (\LE (R_{n}')\subset \ereg )|\leq } &  & \\
 & \qquad  & \leq C\left(\frac{\epsilon }{s}\right)^{-c}+\rho _{3}(\epsilon +s,2^{-n};a,\ereg ,\dreg )+\rho _{3}(\epsilon +s,2^{-n};a,\ereg ,P)
\end{eqnarray*}
where $R_{n}$ is as in lemma \ref{lem:limit} and $R_{n}'$ is a
random walk stopped on $\partial P$. $\rho _{3}(P)$ can be estimated
as in lemma \ref{lem_cont_E_move} to give\[
\rho _{3}(\epsilon +s,2^{-n};a,\ereg ,P)\leq 2\rho _{3}(\epsilon +2s,2^{-n};a,\ereg ,\dreg )+C\left(\frac{d(a,\partial \dreg \cup \partial \ereg )}{s}\right)^{-c}\quad .\]
 This --- with lemma \ref{lem:limit} --- reduces our lemma to an
exercise in calculus.
\end{proof}
\begin{lem}
\label{lem:rho1_OK}For any $\dreg \subset \CC $ open, bounded and
finitely connected, and for any $a\in \dreg $,\[
\lim _{(\epsilon ,\delta )\rightarrow (0,0)}\rho _{1}(\epsilon ,\delta ;a,\dreg )=0\]

\end{lem}
\begin{proof}
Let $\left\{ H_{i}\right\} _{i=1}^{m}$ be the connected components
of $\CC \setminus \dreg $. Let \[
\tau :=\min \{\diam H_{i}\, :\, \diam H_{i}\neq 0\}\quad .\]
A simple discretization estimates shows that for $\epsilon <\tau $
and $\delta $ sufficiently small we can ignore the holes with positive
diameter. For the punctures, let $H$ be the set of punctures and
let $d=d(a,H)$. Then\[
q(a,\partial H,\partial \dreg ,\delta \ZZ ^{2})\leq Cm\frac{\log \diam \dreg -\log d}{\log \delta ^{-1}+\log \diam \dreg }\]
which clearly converge to $0$ as $\delta \rightarrow 0$.
\end{proof}
\begin{defn*}
The family of polygons $\ereg $ such that \begin{equation}
\rho (\epsilon ,2^{-n};a,\ereg ,\dreg )\rightarrow 0\label{rho_sane}\end{equation}
 will be denoted by $\eclass $.
\end{defn*}
Note that $\eclass $ is closed to finite unions and intersections.

\begin{lem}
\label{lem:denseO}$\eclass $ is dense in the sense that for any
bounded open $O$ and for any $\epsilon $ there exists a set $V\in \eclass $
such that $O\subset V\subset O+\epsilon \DD $.
\end{lem}
\begin{proof}
In view of lemma \ref{lem:rho1_OK} we need only estimate $\rho _{2}$.
Assume $\ereg $ is a polygon that satisfies the additional requirement
that each final segment of $\ereg $ before meeting $\partial \dreg $
is a segment from a point $\xi $ to the point of $\partial \dreg $
closest to $\xi $. Then we have, for $r$ sufficiently small, that
$X_{2}$ (from the definition of $\rho _{2}$, (\ref{eq:defX2_forrho2}))
is contained in a finite union of balls of radius $2r$. This obviously
gives that $\rho _{2}\rightarrow 0$. Finally, it is a fun topological
exercise to see that the additional condition above does not interfere
with the density of the family of polygons (in the same sense as above).
\end{proof}
\begin{defn*}
In $\Haus (\bar{\dreg })$ define \[
Y(E_{0};E_{1},...,E_{n}):=\{F\in \Haus (\bar{\dreg })\, :\, F\subset E_{0}\; \wedge \; (F\not \subset E_{i}\; \forall i\geq 1)\}\]
and\[
\ieeclass :=\left\{ \bigcup _{i=1}^{k}Y(E_{0}^{i};E_{1}^{i},...,E_{n_{i}}^{i})\, :\, E_{j}^{i}\in \eclass \right\} \quad .\]

\end{defn*}
\begin{lem}
\label{lem:denseOO}$\ieeclass $ is dense in the sense that for every
open set $O\subset \Haus (\bar{\dreg })$ and for every $\epsilon $
there exists $V\in \ieeclass $ such that \[
O\subset V\subset B(O,\epsilon )\]
 where $B(O,\epsilon ):=\{F\in \Haus (\bar{\dreg })\, :\, d(F,O)<\epsilon \}$.
\end{lem}
The proof follows easily from lemma \ref{lem:denseO} and the existence
of $\epsilon $-nets of finite sets in $\Haus (\bar{\dreg })$ and
we shall omit it. 

\begin{lem}
\label{lem:PYConvrg}With the notations of lemma \ref{lem_limit},
for every $Y\in \ieeclass $ \[
\PP (\LE (R_{n})\in Y)\]
converges.
\end{lem}
\begin{proof}
This follows from lemma \ref{lem_limit}, the inclusion-exclusion
principle and some set algebra.
\end{proof}
\begin{thm*}
\label{the_theorem}Let $\dreg \subset \CC $ be an open bounded finitely
connected set, let $a\in \dreg $ and let $R_{n}$ be a random walk
on $G:=2^{-n}\ZZ ^{2}$ starting from $a$ and stopped on $\partial _{G}\dreg $.
Let $\mu _{n}$ be the distribution measures of $\LE (R_{n})$. Then
$\mu _{n}$ converge in the weak-{*} topology of $\Meas (\Haus (\bar{\dreg }))$.
\end{thm*}
\begin{proof}
It is enough to show that $\mu _{n}(f)$ converges for every continuous
$f$ since this proves that $\mu _{n}$ converges to an arbitrarily
chosen sub-sequence limit. Let $f$ be a continuous function on $\Haus (\bar{\dreg })$.
Since $\Haus (\bar{\dreg })$ is compact $f$ is bounded and uniformly
continuous and we may write, for every $N$, and for $M<\min f$,
\[
f=M+\sum _{k=MN}^{\infty }\frac{1}{N}\one _{O_{k}}+E\]
where $O_{k}:=f^{-1}\left(\left[\frac{k}{N},\infty \right[\right)$,
the sum is finite and $|E|\leq \frac{1}{N}$. Further, the uniform
continuity of $f$ gives that for some $\epsilon $, $B(O_{k+1},\epsilon )\subset O_{k}$
for all $k$. Lemma \ref{lem:denseOO} allows us to take $V_{k}\in \ieeclass $
satisfying \[
O_{k}\subset V_{k}\subset B(O_{k},\epsilon )\]
 and get from lemma \ref{lem:PYConvrg} that $\mu _{n}(M+\sum \one _{V_{k}})$
converges and \[
\left|f-M-\sum \one _{V_{k}}\right|\leq \frac{2}{N}\]
 so\[
\limsup \mu _{n}(f)-\liminf \mu _{n}(f)\leq \frac{2}{N}\]
and since this is true for all $N$ and for all $f$ the theorem is
proved.
\end{proof}

\subsection{\label{chap:extensions}Extensions}

The technique demonstrated in this paper is quite flexible. The only
property of loop-erased random walk crucially used is symmetry. Below
are a few possible future directions.

\begin{itemize}
\item A hybrid graph that interpolates between $\frac{1}{2N}\ZZ ^{2}$ and
$\frac{1}{3N}\ZZ ^{2}$ can be used to show that the scaling limit
is invariant to multiplication by $\frac{2}{3}$ (i.e. if $L(\dreg )$
is the scaling limit on $\dreg $ then $L(\dreg )\sim \frac{2}{3}L(\frac{3}{2}\dreg )$.
This will easily give that the limit of loop-erased random walks on
$\dreg \cap \delta \ZZ ^{2}$ converges to a weak limit as $\delta \rightarrow 0$
continuously.
\item It is possible to use this technique to show that the scaling limit
is invariant to conformal maps. Very roughly, the proof is as follows:
it is only necessary to define a hybrid graph that interpolates between
$\frac{1}{N}\ZZ ^{2}$ and $\varphi \left(\frac{1}{N}\ZZ ^{2}\right)$
where $\varphi $ is the conformal map. If $\varphi $ is close to
1 in the sense that $|\varphi '-1|\leq \epsilon $ and $|\varphi ''|\leq \epsilon $
then it is possible to construct the graph by linking points on the
seams to the closest points on the other part of the graph. The requirement
that $\bar{\bar{G}}$ is only within $O(\frac{1}{N})$ distance from
$\ZZ ^{2}$ gives linear equations for the weights of these links
which can always be solved and the solution is bounded. This reduces
the calculation of $\beta $ (i.e.~the proof of lemma \ref{lemma_a_wired})
to a few specific graphs.
\item I believe this technique might work in 3 dimensions as well. We are
now working on the details.
\item On the other hand, it is hard to image how this technique might be
used for percolation, the UST Peano curve, or any other process where
quasi-loops do exist (in other words, where the limit is $\textrm{SLE}_{\kappa }$
with $\kappa >4$).
\end{itemize}

\appendix

\section{\label{chap_computer}\label{sect:hyrbid_numerology}Proof of (\ref{two_seams_est})}

There is nothing much to say here, really. Clearly we can assume $N=1$.
The values of the harmonic potential of $\mathbb{Z}^{2}$ at specific
points can be calculated by McCrea-Whipple's algorithm. This algorithm
basically uses the fact that there is a close formula for $a(n+in)$,
namely\[
a(n+in)=\frac{1}{\pi }\left(1+\frac{1}{3}+\dotsb +\frac{1}{2n-1}\right)\quad .\]
With these values at hand, the value at any other point can be calculated
using a line-by-line recursion which uses only the harmonicity and
symmetry to $\frac{\pi }{4}$ rotations. See \cite[chapter 15]{S76}
for a more detailed exposition. This allows to calculate $\sum \Delta b_{v}$
on a finite rectangle ($200$ was used for the results below). To
estimate the error outside this rectangle one needs to explicitize
the constants in the proofs of sublemmas \ref{sublem:est_A}-\ref{sublemma_two_seams}.
It must be noted, though, that this requires to know a value for $C_{\ref {a2ndord}}$,
the constant in the estimate (\ref{asimp_a}) of the harmonic potential
on $\mathbb{Z}^{2}$. This is done in \cite[section 4]{KS}, and the
value is \begin{equation}
C_{\ref {a2ndord}}=\frac{9}{4}\left(17-\frac{48+\log 72+2\gamma }{\pi }\right)=0.0172...\label{Val_C1}\end{equation}
The \vpageref[following table][table ]{table:betav} summarizes the
values of the maximal $\beta _{v}:=\sum |\Delta b_{v}-\delta _{v}|$
and the $v$ where they occur for all configurations of $D$ in sublemma
\ref{sublemma_two_seams}. All numerical results are with an error
of $\pm 0.02$.%
\begin{table}
[h]\begin{tabular}{|c|c|c|}
\hline 
$D$&
$\max \beta _{v}$&
happens at\\
\hline
\hline 
$\ZZ ^{2}$ or $\emptyset $&
$0$&
everywhere\\
\hline 
$\ZZ ^{2}\setminus \left(\ZZ ^{-}\right)^{2}$ &
$0.31$&
$-1-i$\\
\hline 
$\ZZ +i\ZZ ^{+}$&
$0.19$&
$\ZZ +\frac{1}{2}-i$\\
\hline 
$\left(\ZZ ^{+}\right)^{2}\cup \left(\ZZ ^{-}\right)^{2}$  &
$0.34$&
$-\frac{3}{2},-\frac{3}{2}i,\frac{1}{2}-i,-1+\frac{1}{2}i$\\
\hline 
$\left(\ZZ ^{+}\right)^{2}$&
$0.39$&
$-\frac{1}{2}-i,-1-\frac{1}{2}i$\\
\hline
\end{tabular}

\caption{\label{table:betav}$\beta $ for simple hybrid graphs}
\end{table}
All programs used are available upon demand.

\end{document}